%% file: art.tex
\documentclass[reqno,11pt]{amsart}

\input mac.tex

\begin{document}  

\nocite{*}

\title[Generalized dendriform algebras]{Polarization and deformations of generalized dendriform algebras}
\author[Ospel]{Cyrille Ospel} \address{Cyrille Ospel, LaSIE, UMR 7356 CNRS, Université de La Rochelle,
  Av. M. Crépeau, 17042 La Rochelle cedex 1, La Rochelle, France}\email{cospel@univ-lr.fr}

\author[Panaite]{Florin Panaite}
\address{Florin Panaite, Institute of Mathematics of the Romanian Academy, PO-Box 1-764, RO-014700 Bucharest, Romania}
\email{florin.panaite@imar.ro}

\author[Vanhaecke]{Pol Vanhaecke}
\address{Pol Vanhaecke, LMA, UMR 7348 CNRS,
  Universit\'e de Poitiers, 11 BD Marie et Pierre Curie, Bât.\ H3 - TSA 61125,
  86 073 Poitiers CEDEX 9, France}\email{pol.vanhaecke@math.univ-poitiers.fr}



\begin{abstract}
  We generalize three results of M. Aguiar, which are valid for Loday's dendriform algebras, to arbitrary
  dendriform algebras, i.e., dendriform algebras associated to algebras satisfying any given set of relations. We
  define these dendriform algebras using a bimodule property and show how the dendriform relations are easily
  determined. An important concept which we use is the notion of polarization of an algebra, which we generalize
  here to (arbitrary) dendriform algebras: it leads to a generalization of two of Aguiar's results, dealing with
  deformations and filtrations of dendriform algebras. We also introduce weak Rota-Baxter operators for arbitrary
  algebras, which lead to the construction of generalized dendriform algebras and to a generalization of Aguiar's
  third result, which provides an interpretation of the natural relation between infinitesimal bialgebras and
  pre-Lie algebras in terms of dendriform algebras. Throughout the text, we give many examples and show how they
  are related.
\end{abstract}

\maketitle

\setcounter{tocdepth}{1}

\tableofcontents

\input intro.tex
\input dendri.tex
\input rb.tex

\input pol.tex

\bibliographystyle{abbrv} \bibliography{ref}

\end{document}

%% file: mac.tex
\usepackage{graphicx}         
\usepackage{amsfonts}
\usepackage{scalerel,amssymb}
\usepackage[latin1]{inputenc}
\usepackage[all]{xy}
\usepackage{tikz-cd}        
\usepackage{pgf,tikz}

\newtheorem{thm}{Theorem}[section]
\newtheorem{cor}[thm]{Corollary}
\newtheorem{lemma}[thm]{Lemma}

\newtheorem{prop}[thm]{Proposition}
\theoremstyle{definition}
\newtheorem{defn}[thm]{Definition}
\theoremstyle{remark}

\newtheorem{remark}[thm]{Remark}
\newtheorem{example}[thm]{Example}

\renewcommand{\l}{\lambda}
\newcommand{\s}{\sigma}

\newcommand{\PB}{\left\{\cdot\,,\cdot\right\}}

\newcommand{\po}[2]{#1\circ #2}

\newcommand{\lb}[1]{\[#1\]}
\newcommand{\LB}{[\cdot\,,\cdot]}

\newcommand{\gb}[1]{[\![#1]\!]}
\newcommand{\GB}{[\![\cdot\,,\cdot]\!]}

\renewcommand{\(}{\left(}
\renewcommand{\)}{\right)}
\renewcommand{\[}{\left[}
\renewcommand{\]}{\right]}

\newcommand{\set}[1]{\left\{#1\right\}}

\newcommand{\cat}[3]{#1^{\hbox{\tiny #2}}_{\hbox{\tiny #3}}}
\newcommand{\catdend}{\cat{\cC}{dend}{}}
\newcommand{\cattrid}{\cat{\cC}{trid}{}}
\newcommand{\catpol}{\cat{\cC}{}{pol}}
\newcommand{\ucatpol}{\cat{\uC}{}{pol}}
\newcommand{\catdendpol}{\cat{\cC}{dend}{pol}}
\newcommand{\ucatdendpol}{\cat{\uC}{dend}{pol}}
\newcommand{\cattridendpol}{\cat{\cC}{trid}{pol}}

\newcommand{\catpdend}{\cat{\cC}{$'$dend}{}}
\newcommand{\catp}{\cC'}
\newcommand{\catcom}{\cat{\cC}{}{com}}
\newcommand{\catdendcom}{\cat{\cC}{dend}{com}}

\newcommand{\catPdend}{\cat{\cP}{dend}{}}
\newcommand{\catPtrid}{\cat{\cP}{trid}{}}
\newcommand{\catPpol}{\cat{\cP}{}{pol}}
\newcommand{\catPdendpol}{\cat{\cP}{dend}{pol}}
\newcommand{\catPtridendpol}{\cat{\cP}{trid}{pol}}

\newcommand{\cC}{\mathcal C}
\newcommand{\uC}{\underline{\cC}}

\newcommand{\cM}{\mathcal M}
\newcommand{\cP}{\mathcal P}

\newcommand{\cR}{\mathcal R}
\newcommand{\uR}{\underline\cR}

\newcommand{\bbN}{\mathbb N}

\newcommand{\tcM}{\tilde\cM}
\newcommand{\tM}{\tilde M}

\def\mcirc{\mathbin{\scalerel*{\circ}{j}}}

\newcommand{\Ah}{A^\h}

\newcommand{\Rh}{R^\h}
\newcommand{\h}{\nu}
\newcommand{\squ}{\,\square\;}
\newcommand{\fR}{\mathfrak R}
\newcommand{\fS}{\mathfrak S}
\newcommand{\ul}[2]{\underline{#1}_{#2}}
\renewcommand{\a}[1]{\ul a{#1}}
\renewcommand{\b}[1]{\ul b{#1}}
\renewcommand{\c}[1]{\ul c{#1}}
\newcommand{\RF}{R}

\newcommand{\Gr}{\mathop{\rm gr}\nolimits}

\newcommand{\leqs}{\leqslant}
\newcommand{\geqs}{\geqslant}

\newcommand{\Id}{\mathrm{Id}}

\DeclareMathOperator{\AYB}{AYB}
\allowdisplaybreaks[4]

\newif\ifprivate
\privatefalse

 \numberwithin{equation}{section}

\def\???{\ifprivate {\bf {???}} \marginpar{{\Huge {\bf ?}}}\else \fi}
\numberwithin{equation}{section}

%% file: intro.tex
\section{Introduction}
Dendriform algebras were introduced by J.-L. Loday~in \cite{lodaydend} as a dichotomized version of associative
algebras. By definition, a \emph{Loday dendriform algebra} is an algebra $(A,\prec,\succ)$ satisfying, for all
$a,b,c\in A$, the~relations%
\begin{align*}
  (a\prec b)\prec c&=a\prec (b\prec c+b\succ c)\;, \\
  (a \succ b)\prec c&=a\succ (b\prec c)\;,  \\
  (a\prec b+a\succ b)\succ c&=a\succ (b\succ c)\;.
\end{align*}
Summing up these three equations and setting $a\star b:=a\prec b+a\succ b$ for all $a,b\in A$, one sees that
$\star$ is associative, so that $(A,\star)$ is an associative algebra, on which the dendriform operations provide
some extra structure.  In \cite{aguiar}, M. Aguiar introduced the notion of deformation for a
\emph{commutative} dendriform algebra $(A,\prec,\succ)$, where commutativity means that $a\prec b=b\succ a$ for all
$a,b\in A$. He shows that such a deformation makes $(A,\times,\mcirc)$ into a pre-Poisson algebra, where $\times$
stands for $\succ$ and where the new product $\mcirc$ on $A$ is constructed from the first order deformation terms
of $\prec$ and $\succ$.  The notion of a pre-Poisson algebra was also introduced in loc.\ cit.: $(A,\times,\mcirc)$
is a \emph{pre-Poisson} algebra if the following relations are satisfied, for all $a,b,c\in A$:
\begin{eqnarray*}
  a\times (b\times c)&=&(a\times b+b\times a)\times c\;,\\
  (a\times b+b\times a)\mcirc c&=& a\times (b\mcirc c)+b\times (a\mcirc c)\;,\\
  (a\mcirc b-b\mcirc a)\times c&=&a\times (b\mcirc c)-b\mcirc(a\times c)\;, \\
  (a\mcirc b-b\mcirc a)\mcirc c&=&a\mcirc(b\mcirc c)-b\mcirc(a\mcirc c)\;.
\end{eqnarray*}%
This result is a dendriform version of the well-known fact that the skew-symmetrization of the first deformation
term of a deformation of an associative algebra $A$ is a Poisson bracket on $A$. Aguiar also establishes a similar
result for filtered dendriform algebras, also a dendriform version of a well-known result.  Even if these
results can easily be proven by a direct computation, these computations lack a conceptual understanding, which we
will provide in this paper. We do this by generalizing these results to arbitrary dendriform algebras; a key
element is the notion of polarization for (arbitrary) dendriform algebras, which we will introduce.


We define generalized dendriform algebras as follows. Let $\cC$ denote the category of all algebras $(A,\mu)$ which
satisfy a given set of relations $\cR_1=0,$ $\dots,\cR_k=0$. An algebra $(A,\prec,\succ)$ is said to be a
\emph{$\cC$-dendriform algebra} if $(A\times A,\boxtimes)\in\cC$, where $\boxtimes$ is defined for $(a,x),(b,y)\in
A\times A$, by $(a,x)\boxtimes (b,y):=(a\prec b+a\succ b,a\succ y+x\prec b)$. This property can also be expressed
as a bimodule property. The $\cC$-dendriform algebras form a category $\catdend$ with algebra homomorphisms as
morphisms. Taking as relation associativity, we recover the definition of a Loday dendriform algebra. We show that
when all relations are multilinear, the relations which every $\cC$-dendriform algebra must satisfy are easily
obtained from the relations $\cR_i=0$. Generalized dendriform algebras have already been considered from the
operadic point of view in \cite{BaiBellierGuoNi}, but we will not use or need this formalism since the phenomena
and properties which we present are most naturally expressed in terms of the basic algebraic language which we
use.

In order to construct (interesting) examples of generalized dendriform algebras, we introduce the notion of a weak
Rota-Baxter operator. Given any algebra $A$ (satisfying a given set of relations), a linear map $\fR:A\to A$ is
said to be a \emph{weak Rota-Baxter operator} of $A$ if, for all $a,b\in A$, the element
\begin{equation}\label{eq:weak_rb_intro}
    \fR(a\fR(b)+\fR(a)b)-\fR(a)\fR(b)
  \end{equation}%
commutes with all elements of $A$. We show that $(A,\prec,\succ)$ becomes a generalized dendriform algebra upon
setting $a\succ b:=\fR(a)b$ and $a\prec b:=a\fR(b)$ for all $a,b\in A$. More precisely, we show which relations
$\prec$ and $\succ$ will satisfy. When (\ref{eq:weak_rb_intro}) is zero for all $a,b\in A$ (in which case $\fR$ is
called a \emph{Rota-Baxter operator}) these relations are precisely the dendrification of the relations satisfied
by $A$; the same is true for arbitrary weak Rota-Baxter operators in case the relations can be written in
commutator form (see Section \ref{par:dendri_from_weak_RB} for the definition of this notion). As an application,
we generalize yet another result by M. Aguiar \cite{aguiar_inf}, which states that the natural functor which
associates to any $\epsilon$-bialgebra $(A,\mu,\Delta)$ the corresponding pre-Lie algebra $(A,\circ)$, restricted
to the category of quasi-triangular $\epsilon$-bialgebras, admits a natural factorization through the category of
dendriform algebras, i.e., the following diagram is commutative (see Section \ref{par:coboundary} for details):
\begin{center}
  \begin{tikzcd}[row sep=8ex, column sep=15ex]
    \hbox{QT }\epsilon\hbox{-bialg},\mu,r\arrow{d}{\sum_iau_ibv_i,\sum_iu_iav_ib}\arrow{r}{r\cdot a-a\cdot r}
      &\epsilon\hbox{-bialg},\mu,\Delta\arrow{d}{\sum_{(b)}b_{(1)}ab_{(2)}}\\
    \hbox{Assoc}^{\hbox{\tiny dend}},\prec,\succ\arrow{r}{a\succ b-b\prec a}&\hbox{pre-Lie},\circ
  \end{tikzcd}
\end{center}
We show that this diagram can be generalized to coboundary $\epsilon$-bialgebras by replacing in it the two
leftmost entries by coboundary $\epsilon$-bialgebras and $A_3$-dendriform algebras, without changing the arrows;
recall from \cite{gozeremm} that an algebra $(A,\mu)$ is said to be $A_3$-associative if for all $a,b,c\in A$,
\begin{equation*}
  (ab)c+(bc)a+(ca)b=a(bc)+b(ca)+c(ab)\;;
\end{equation*}%
the corresponding dendriform algebras are called \emph{$A_3$-dendriform algebras}.

For algebras $(A,\mu)$ (with one operation), the notion of polarization has been introduced in \cite{MarklRemm}:
the product $\mu$ is decomposed into its symmetric and antisymmetric parts, yielding an algebra $(A,\cdot,\LB)$ for
which the relations are obtained from the relations satisfied by $\mu$. This definition is easily adapted to
generalized dendriform algebras, as indicated in the following commutative diagram of categories and functors,
where the horizontal arrows are isomorphisms of categories (see Section \ref{par:gendend2} for the notation):
\begin{equation*}
  \begin{tikzcd}[row sep=8ex, column sep=15ex]
    \cC,\mu\arrow[shift left]{r}{(ab+ba)/2,(ab-ba)/2}&
      \catpol,\cdot,\LB\arrow[shift left]{l}{a\cdot b+\lb{a,b}}\\
      \catdend,\prec,\succ\arrow{u}{a\prec b+a\succ b}
                          \arrow[shift left]{r}{\frac{a\succ b+b\prec a}2,\frac{a\succ b-b\prec a}2}
    &\catdendpol,*,\circ\arrow[swap]{u}{a*b+b*a,a\circ b-b\circ a}
    \arrow[shift left]{l}{b*a-b\circ a,a*b+a\circ b}
  \end{tikzcd}
\end{equation*}

Thanks to the commutativity of the diagram we can obtain the relations of $\catdendpol$ from the relations of
$\catpol$ in case the latter are multilinear. For example, if $\catpol$ denotes the category of Poisson algebras,
$\catdendpol$ is the category of pre-Poisson algebras, mentioned above.

We give two applications of polarization. Suppose that $(A[[\h]],\prec,\succ)$ is a deformation of a
commutative algebra $(A,\prec_0,\succ_0)$ in $\catdend$ and consider the algebra $(A,\times,\mcirc)$, where
$\times$ stands for $\succ_0$ and $\mcirc$ is defined for $a,b\in A$ by
\begin{equation*}
    a\mcirc b:=\frac{a\succ b-b\prec a}{2\h}{\Big\vert_{\h=0}}\;.
\end{equation*}
We show that $(A,\times,\mcirc)\in\ucatdendpol$, where $\ucatdendpol$ is the category of all polarized dendriform
algebras $(A,*,\circ)$ satisfying for all relations $\cR=0$ of $\catdendpol$ the relation $\uR=0$; here, $\uR$
stands for the lowest weight part of $\cR$, where the weight of a monomial in $A$ is defined as being the number of
operations $\circ$ that it contains. Notice that these relations can easily be computed. A prime example is the
case in which $(A,\prec,\succ)$ is a Loday dendriform algebra; we then recover Aguiar's result since then
$\ucatdendpol$ is the category of pre-Poisson algebras. As a second application, we show that a similar result
holds for filtered commutative algebras in $\catdend$. Both applications admit also an anticommutative version.

The structure of the paper is as follows. We introduce in Section \ref{sec:dendri} the notion of a $\cC$-dendriform
algebra and we show how the relations satisfied by all $\cC$-dendriform algebras can be obtained from the relations
in $\cC$. We give several examples and we show how they are related, both in their original and in their dendrified
form. Rota-Baxter and weak Rota-Baxter operators are shown in Section~\ref{sec:A3} to provide constructions of
$\cC$-dendriform algebras; we give an application of it to $\epsilon$-bialgebras. The notion of polarization for
dendriform algebras is introduced in Section~\ref{par:gendend2} and again we show how for polarized algebras,
defined by multilinear relations, the relations satisfied by the corresponding dendriform algebras are obtained. As
an application, we give a conceptual proof of the generalization to $\cC$-dendriform algebras of Aguiar's results,
stated at the beginning of the introduction; this yields, in particular, a conceptual proof of these results. All
our results extend to $\cC$-tridendriform algebras; throughout the paper, we will indicate these generalizations in
some short remarks.

\noindent{\bf Conventions.} All algebraic structures which we consider (algebras, modules, bialgebras, etc.)  are
defined over a commutative ring $\RF$ in which~2 is invertible. Since the base ring $R$ will never change, we
denote $\otimes_{\RF}$ simply by~$\otimes $.  By ``$R$-algebra'', which we call simply ``algebra'', we mean an
$(n+1)$-tuple $(A, \mu_1,\dots,\mu_n)$, where $A$ is an $R$-module and $\mu_i :A\otimes A\rightarrow A$ is a
\emph{product}, i.e., a linear map, for $i=1,\dots,n$.  By an \emph{algebra homomorphism} between two algebras $(A,
\mu_1,\dots,\mu_n)$ and $(A', \mu_1',\dots,\mu_n')$ we mean a linear map $f:A\to A'$ such that $f(\mu_i(a\otimes
b))=\mu_i'(f(a)\otimes f(b))$ for all $a,b\in A$ and all $i=1,\dots,n$. Unless otherwise specified, the products
$\mu_i$ are not assumed to have any extra properties.
%
%
In the case of an algebra $(A,\mu)$ with one product, we usually write $ab$ for $\mu (a\otimes b)$.
We use the standard notations $S_n$ and $A_n$ for the symmetric and alternating groups of degree $n$.

%% file: dendri.tex
\section{Dendriform algebras}\label{sec:dendri}
In this section, we recall the notion of a Loday dendriform algebra and we show that it naturally generalizes to
algebras $(A,\mu)$, defined by any finite collection of relations; Loday dendriform algebras correspond to
associative algebras, which are defined by a single relation, namely associativity.  We also show that when the
relations in the original algebra are multilinear, the relations which hold in the corresponding dendriform
algebras are easily determined. Recall that we write $ab$ for~$\mu(a\otimes b)$.

\subsection{Loday dendriform algebras}\label{par:classical_dendri}
We first recall from \cite{lodaydend} the notion of a Loday dendriform algebra.
\begin{defn}\label{def:dendri}
A \emph{Loday dendriform algebra} is an algebra $(A,\prec,\succ)$ satisfying for all $a, b, c\in A$ the following
relations:
\begin{align}
  (a\prec b)\prec c&=a\prec (b\prec c+b\succ c)\;,  \label{dend1} \\
  (a \succ b)\prec c&=a\succ (b\prec c)\;, \label{dend2} \\
  (a\prec b+a\succ b)\succ c&=a\succ (b\succ c)\;. \label{dend3}
\end{align}
\end{defn}
The terminology \emph{dendriform} comes from the shape of the free Loday dendriform algebra, which is naturally
described in terms of planar binary trees [loc.\ cit.\ Sections 5.4 and 5.7]. Dendriform algebras can be considered
as a dichotomized version of an associative algebra: defining $a\star b:=a\prec b+a\succ b$ for all $a,b\in A$, the
newly formed algebra $(A,\star)$ is associative. In fact, Loday dendriform algebras can be characterized as follows
(see \cite{Bai_Guo_Ni}):
\begin{prop}\label{prp:dendri_char}
  Let $(A,\prec,\succ)$ be an algebra and let $\star$ denote the sum of $\prec$ and $\succ$. Then $(A,\prec,\succ)$
  is a dendriform algebra if and only if the following conditions are satisfied:
  \begin{enumerate}
  \item[(1)] $(A,\star)$ is an associative algebra;
  \item[(2)] $(A,\succ,\prec)$ is an $(A,\star)$-bimodule.
  \end{enumerate}
\end{prop}
In this characterization, the notion of bimodule (over an associative algebra) is the standard one; see the lines
following Definition \ref{def:dendri_gen} below for the more general concept of a bimodule over other types of
algebras.

Conditions (1) and (2) of the proposition can be restated as the single condition that $(A\times A,\boxtimes)$ is
associative, where the product $\boxtimes$ is defined, for $(a,x),(b,y)\in A\times A$, by
\begin{equation}\label{eq:box_def}
  (a,x)\boxtimes (b,y):=(a\star b,a\succ y+x\prec b)\;.
\end{equation}%
The proof of the equivalence is a direct consequence of the following formulas, valid for all $(a,x),(b,y),(c,z)\in
A\times A$:
\begin{align*}
  ((a,x)\boxtimes(b,y))\boxtimes(c,z)&=((a\star b)\star c,(a\star b)\succ z+(a\succ y)\prec c+(x\prec b)\prec c)\;,\\
  (a,x)\boxtimes((b,y)\boxtimes(c,z))&=(a\star(b\star c),a\succ(b\succ z)+a\succ(y\prec c)+x\prec(b\star c))\;.
\end{align*}
It follows that a Loday dendriform algebra can equivalently be defined as an algebra $(A,\prec,\succ)$ such that
$(A\times A,\boxtimes)$ is associative, where $\boxtimes$ is defined by~(\ref{eq:box_def}). It is this more
conceptual definition which we will generalize.

\subsection{$\cC$-dendriform algebras}\label{sec:C-dendri}
Let $\cR_1=0,\dots,\cR_k=0$ be given relations and denote by $\cC$ the category of all algebras (with one
operation) which satisfy these relations called, by a slight abuse of language, \emph{the relations
of~$\cC$}. Morphisms in $\cC$ are algebra homomorphisms. If $(A,\mu)$ is an object of~$\cC$, we write
$(A,\mu)\in\cC$.

\begin{defn}\label{def:dendri_gen}
  An algebra $(A,\prec,\succ)$ is said to be a \emph{$\cC$-dendriform algebra} if $(A\times A,\boxtimes)\in\cC$,
  where $\boxtimes$ is defined for $(a,x),(b,y)\in A\times A$, by
  \begin{equation}\label{eq:semi-direct}
    (a,x)\boxtimes (b,y):=(a\prec b+a\succ b,a\succ y+x\prec b)\;.
  \end{equation}%
\end{defn}
Taking $x=y=0$ in (\ref{eq:semi-direct}), it is clear that if $(A,\prec,\succ)$ is a $\cC$-dendriform algebra, then
$(A,\star)\in\cC$, where $\star$ denotes the sum of $\prec$ and $\succ$.  In the language of bimodules (for general
algebras, not necessarily associative), the property that $(A\times A,\boxtimes)$ belongs to $\cC$, where
$\boxtimes$ is defined by (\ref{eq:semi-direct}), is \emph{by definition} precisely the condition that
$(A,\star)\in\cC$ and that $(A,\succ,\prec)$ is an $(A,\star)$-bimodule with respect to $\cC$ (see
\cite{schafer}). 
\begin{remark}\label{rem:tridendri}
Definition \ref{def:dendri_gen} admits the following natural generalization: using the notations and under the
assumptions of that definition, an algebra $(A,\prec,\succ,.)$ is said to be a \emph{$\cC$-tridendriform algebra}
if $(A\times A,\boxtimes)\in\cC$, where $\boxtimes$ is now defined for $(a,x),(b,y)\in A\times A$, by
  \begin{equation}\label{eq:semi-direct_tri}
    (a,x)\boxtimes (b,y):=(a\prec b+a\succ b+a.b,a\succ y+x\prec b+x.y)\;.
  \end{equation}
In the particular case when $a.b=0$ for all $a,b\in A$, one recovers the above definition of a $\cC$-dendriform
algebra. Also, taking for $\cC$ the category of all associative algebras, one recovers the classical notion of a
tridendriform algebra, as first introduced by J.-L. Loday and M. Ronco in \cite{lodayronco} (for a proof, see
\cite{Bai_Guo_Ni} in which our definition of a $\cC$-tridendriform algebra appears in the associative case as a
characterization of a tridendriform algebra).

\end{remark}

\subsection{Algebras defined by multilinear relations}\label{sec:multilinear}
The relations, satisfied by the algebras which we will consider, are multilinear and we will show how for such
relations we can easily obtain the corresponding relations which must be satisfied by the corresponding dendriform
algebras; we do this for one relation at a time.  Our method is based on the fact that, by multilinearity, the
condition that $(A\times A,\boxtimes)$ belongs to $\cC$ is equivalent to the conditions obtained by demanding that
the relations are satisfied for all possible $n$-tuplets (for an $n$-linear relation) of elements of $A\times A$,
taken from a generating set of $A\times A$. We take this generating set to be the union of $A_0:=A\times\set0$ and
$A_1:=\set0\times A$.  We will find it convenient to use for any $a\in A$ the following notation: $\a0:=(a,0)$ and
$\a1:=(0,a)$; also, when we consider elements $\a0\in A_0$ or $\a1\in A_1$ we implicitly assume that $a\in A$. In
this notation, (\ref{eq:semi-direct}) is equivalenty described by the following table, in which $a$ and $b$ stand
for arbitrary elements of $A$:
\begin{table}[h]
  \def\arraystretch{1.8}
  \setlength\tabcolsep{0.7cm}
  \centering
\begin{tabular}{c|cc}
   $\boxtimes$&$\b0$&$\b1$\\
  \hline
  $\a0$&$\ul{a\star b}0$&$\ul{a\succ b}1$\\
  $\a1$&$\ul{a\prec b}1$&$(0,0)$\\
\end{tabular}
\medskip
\caption{The product $\boxtimes$ for generators of $A\times A$.}\label{tab:boxtimes}
\end{table}
%
%

\vskip -0.5cm
We explain the procedure in the case of a trilinear relation, the case of a bilinear relation being too
simple\footnote{When $R$ is a field, the only non-trivial bilinear relations are commutativity and
anticommutativity.}  to illustrate how it works; see 
Remark~\ref{rem:n-linear} below for the case of an $n$-linear relation. By a \emph{trilinear relation} on an
algebra $(A,\mu)$, we mean a relation in three variables which is linear in each of the variables, i.e., the
relation is of the form $\cR=0$, where
\begin{equation}\label{eq:trilinear}
  \cR(a_1,a_2,a_3)=\sum_{\s\in S_3}\l_\s\;(a_{\s(1)}a_{\s(2)})a_{\s(3)}+
                   \sum_{\s\in S_3}\l'_\s\; a_{\s(1)}(a_{\s(2)}a_{\s(3)})\;.
\end{equation}%
The 12 constants $\l_\s$ and $\l_\s'$ belong to the base ring $\RF$. The associativity relation, $(ab)c-a(bc)=0$, is
an example; there are many other such relations, such as the ones defining Leibniz algebras, NAP algebras, pre-Lie
algebras, Lie-admissible algebras, and so on. Several of these, and some others, will be considered below, where
their definition will be recalled.

Let $\cR=0$ be a trilinear relation and let us denote by $\cR_\boxtimes$ (resp.\ $\cR_\star$) the formula $\cR$ in
which the product $\mu$ is replaced by $\boxtimes$ (resp.\ $\star$). We show how to obtain the corresponding
relations for a $\cC$-dendriform algebra.
%

\noindent$\bullet$\quad
If we take three arbitrary elements $\a0,\b0,\c0$ in $A_0$, then 
\begin{equation*}
  \qquad\quad(\a0\boxtimes \b0)\boxtimes \c0=\ul{(a\star b)\star c}0\;,\ \hbox{and}\
  \a0\boxtimes(\b0\boxtimes \c0)=\ul{a\star(b\star c)}0\;,
\end{equation*}
so that 
%
  $\cR_\boxtimes(\a0,\b0,\c0)=\ul{\cR_\star(a,b,c)}0\;,$
%
for all $a,b,c\in A$. Therefore, the relation which we find is that $\cR_\star=0$, i.e., that $(A,\star)\in\cC$. As
we will see in the next item, this relation needs not be stated explicitly, because it follows from the other
relations.

\noindent$\bullet$\quad When we take two elements in $A_0$ and one in $A_1$, we get from $\cR_\boxtimes=0$ three
non-trivial relations which may be linearly dependent. Notice that
\begin{align*}
  \qquad\lefteqn{(\a0\boxtimes \b0)\boxtimes \c1+(\a0\boxtimes \b1)\boxtimes \c0+(\a1\boxtimes \b0)\boxtimes \c0}\\
      &=\ul{(a\star b)\succ c+(a\succ b)\prec c+(a\prec b)\prec c}1=\ul{(a\star b)\star c}1\;,
\end{align*}
for any $a,b,c\in A$, and similarly with the opposite parenthesizing,
\begin{equation*}
  \quad\qquad\a0\boxtimes(\b0\boxtimes \c1)+\a0\boxtimes(\b1\boxtimes \c0)+\a1\boxtimes(\b0\boxtimes \c0)=
    \ul{a\star(b\star c)}1\;.
\end{equation*}
If we write $\cR$ as in (\ref{eq:trilinear}), then it follows from these two equations that 
\begin{eqnarray}\label{eq:sum_of_all}
  \lefteqn{\cR_\boxtimes(\ul{a_1}0,\ul{a_2}0,\ul{a_3}1)+\cR_\boxtimes(\ul{a_1}0,\ul{a_2}1,\ul{a_3}0)+
    \cR_\boxtimes(\ul{a_1}1,\ul{a_2}0,\ul{a_3}0)}\nonumber\\
  &=&\sum_{\s\in S_3}\l_\s\ul{(a_{\s(1)}\star a_{\s(2)})\star a_{\s(3)}}1+
    \sum_{\s\in S_3}\l_\s\ul{(a_{\s(1)}\star a_{\s(2)})\star a_{\s(3)}}1\nonumber\\
  &&+\sum_{\s\in S_3}\l'_\s\ul{a_{\s(1)}\star (a_{\s(2)}\star a_{\s(3)})}1+
     \sum_{\s\in S_3}\l'_\s\ul{a_{\s(1)}\star (a_{\s(2)}\star a_{\s(3)})}1\nonumber\\
  &=&\ul{\cR_\star(a_1,a_2,a_3)}1\;,
\end{eqnarray}
and so the sum of the three relations which we just found for $\prec$ and $\succ$ is precisely the corresponding
relation for their sum $\star$, as stated above.

\noindent$\bullet$\quad Taking at most one element in $A_0$ and the other ones in $A_1$ gives trivial relations,
because a triple product in $(A\times A,\boxtimes)$ vanishes as soon as at least two of its factors belong to
$A_1$, as follows at once from the definition of $\boxtimes$.

The upshot is that a trilinear relation $\cR=0$ gives rise to at most three independent relations, which are found
by considering $\cR_\boxtimes$ for a triplet of elements in $A\times A$, where two of them are arbitrary elements
in $A_0$ and the other one in $A_1$. Notice that, when~$\cR$ is invariant under a cyclic permutation in its three
variables, the three obtained relations will be the same, so that only one such triplet has to be considered;
similarly, when~$\cR$ is invariant under a transposition of two of the three variables, only two triplets need to
be considered. Since the defining relations of many types of algebras are quite symmetric, we will see below
several examples of this.


\begin{remark}\label{rem:n-linear}
The above analysis is also valid for $n$-linear relations, with $n>3$: in order to obtain all $\cC$-dendriform
relations, it suffices to substitute $n-1$ elements from $A_0$ and one from $A_1$, and this in the $n$ possible
ways. To see this, notice first that if one substitutes in any monomial $a_1\boxtimes a_2\boxtimes\cdots\boxtimes
a_n$ (with any parenthesizing)
at least two elements from $A_1$ and the other ones from $A_0$, one always gets zero, because $A_0\boxtimes A_1$
and $A_1\boxtimes A_0$ are contained in $A_1$ and $A_1\boxtimes A_1=\set{(0,0)}$.  It remains to be shown that the
relation, which is obtained by substituting~$n$ elements from~$A_0$, follows from the $n$ relations which are
obtained by substituting $n-1$ elements from $A_0$ and one element from $A_1$.  Consider a monomial $a_1a_2\dots
a_n$ in $A$, with some parenthesizing, and denote for $i=1,2,\dots,n$,
\begin{align*}
  X&:=\ul{a_1}0\boxtimes\ul{a_2}0\boxtimes\cdots\boxtimes\ul{a_n}0=\ul{a_1\star a_2\star\cdots \star a_n}0\;,\\
  X_i&:=\ul{a_1}0\boxtimes\ul{a_2}0\boxtimes\cdots\boxtimes\ul{a_{i-1}}0\boxtimes\ul{a_{i}}1
    \boxtimes\ul{a_{i+1}}0\boxtimes\cdots\boxtimes \ul{a_n}0\;,
\end{align*}
with the same parenthesizing. Notice that $X\in A_0$ and that $X_i\in A_1$ for $i=1,2,\dots,n.$ Defining $a\in A$
by $X=\ul a0$ (i.e., $a=a_1\star a_2\star\cdots \star a_n$, with the same parenthesizing), we show that
$\sum_{i=1}^n X_i=\ul a1$. We do this by induction on $n$, the case of $n=3$ already being proven above. We can
write~$X$ (uniquely, as dictated by the parenthesizing) as $X=X'\boxtimes X''$, where
\begin{equation*}
  X'=\ul{a_1}0\boxtimes\ul{a_2}0\boxtimes\cdots\boxtimes\ul{a_m}0\;, \qquad
  X''=\ul{a_{m+1}}0\boxtimes\ul{a_{m+2}}0\boxtimes\cdots\boxtimes\ul{a_n}0\;,
\end{equation*}%
with $1\leqslant m<n$, and both $X'$ and $X''$ come with a parenthesizing inherited from the one of~$X$. We define
for $i=1,\dots,m$ (resp.\ for $i=m+1,\dots,n$) the element $X'_i$ (resp.\ $X''_i$) analogously to the definition of
$X_i$ above. If we apply the induction hypothesis to $X'$ and $X''$, we get $\sum_{i=1}^mX_i'=\ul{a'}1$ and
$\sum_{i=m+1}^nX_i''=\ul{a''}1,$ {where} $X'=\ul{a'}0$ and $X''=\ul{a''}0.$
%
It follows that
\begin{align*}  
  \sum_{i=1}^n X_i&=\sum_{i=1}^m X'_i\;\boxtimes X''+X'\boxtimes\sum_{i=m+1}^nX''_i
    =\ul{a'}1\boxtimes\ul{a''}0+\ul{a'}0\boxtimes\ul{a''}1\\
    &=\ul{a'\prec a''}1+\ul{a'\succ a''}1=\ul{a'\star a''}1\;,
\end{align*}
while $X=X'\boxtimes X''=\ul{a'}0\boxtimes\ul{a''}0=\ul{a'\star a''}0,$ so that $\sum_{i=1}^n X_i=\ul a1$ where
$X=\ul a0$. It proves the announced property for $n$-linear relations, for all~$n$.
\end{remark}

\begin{remark}
For relations which are sums of $k$-linear relations, with $k$ varying from $1$ to $n$, the above procedure can be
adapted, but there is no need to do this since for $k=1,\dots,n$ the $k$-linear part of such a relation $\cR=0$ is
itself a relation.  To show this, one shows that the leading ($n$-linear) part is a relation, which follows by
substituting successively $a_i=0$ for $i=1,\dots,n$.
\end{remark}

\begin{remark}\label{rem:tridendri_rel}
For $\cC$-tridendriform algebras (see Remark \ref{rem:tridendri}), where $\cC$ is defined by multilinear relations,
the relations are obtained in the same way as in the case of $\cC$-dendriform algebras, but there will be many more
relations. Indeed, given an $n$-linear relation $\cR=0$, substituting in $\cR_\boxtimes=0$ two or more elements
from $A_1$ and the other ones from $A_0$ will lead to a non-trivial relation, contrary to what we have seen in the
case of a $\cC$-dendriform algebra. We will therefore get $2^n$ relations for a $\cC$-tridendriform algebra, rather
than $n$. It can be shown that the relation, obtained by substituting in $\cR_\boxtimes$ only elements from $A_0$,
is the sum of all $2^n-1$ relations obtained by substituting in~$\cR_\boxtimes$ at least one element from $A_1$ and
the other elements from $A_0$. However, apart from this, these $2^n$ relations are in general independent.
\end{remark}

\subsection{Examples}\label{par:examples}
We illustrate the above procedure in the following examples.
In Section \ref{par:relating_examples}, we will
show how these examples are related.
\begin{example}\label{exa:classical}
We start with the case of a Loday dendriform algebra, which we recalled in Section~\ref{par:classical_dendri}:
here, the only relation is associativity. We show how we obtain the relations of Definition \ref{def:dendri} from
the associativity of $\boxtimes$. First, take $\a0,\b0$ in $A_0$ and $\c1$ in $A_1$. Then, by the associativity of
$\boxtimes$ and by Table \ref{tab:boxtimes},
\begin{equation*}
  \ul{(a\star b)\succ c}1=(\a0\boxtimes \b0)\boxtimes \c1=\a0\boxtimes (\b0\boxtimes \c1)=\ul{a\succ(b\succ c)}1\;,
\end{equation*}%
so that $(a\star b)\succ c=a\succ(b\succ c)$, which is (\ref{dend3}). Relations (\ref{dend2}) and~(\ref{dend1}) are
similarly obtained by taking $\a0,\c0$ in $A_0$ and $\b1$ in $A_1$ (resp.\ $\b0,\c0$ in $A_0$ and $\a1$ in $A_1$).
\end{example}
\begin{example}\label{ex:L-dendri}
A \emph{pre-Lie algebra $(A,\mu)$} is an algebra for which the \emph{associator}, defined by $(a,b,c):=(ab)c-a(bc)$
is symmetric in its first two variables, $(a,b,c)=(b,a,c)$ for all $a,b,c\in A$. Thus, the trilinear relation which
defines pre-Lie algebras is given by
\begin{equation}\label{eq:pre-Lie}
  (ab)c-a(bc)=(ba)c-b(ac)\;.
\end{equation}%
%
Let $\cC_{pL}$ denote the category of all pre-Lie algebras. Using the above procedure, we obtain the relations
which any $\cC_{pL}$-dendriform algebra $(A,\prec,\succ)$ must satisfy, by substituting in the relation
\begin{eqnarray}\label{eq:L_dend_sat}
  \lefteqn{((a,x)\boxtimes (b,y))\boxtimes (c,z)-(a,x)\boxtimes((b,y)\boxtimes(c,z))}\nonumber\\
  &&=((b,y)\boxtimes (a,x))\boxtimes (c,z)-(b,y)\boxtimes((a,x)\boxtimes (c,z))\;,
\end{eqnarray}%
two elements from $A_0$ and one from $A_1$. Substituting $\a0,\,\b0$ and $\c1$ in (\ref{eq:L_dend_sat}), we get,
using Table \ref{tab:boxtimes},
\begin{equation*}
  \ul{(a\star b)\succ c}1-\ul{a\succ(b\succ c)}1=  \ul{(b\star a)\succ c}1-\ul{b\succ(a\succ c)}1\;,
\end{equation*}%
which leads to the relation
\begin{equation}\label{eq:L-dend1}
  (a\star b)\succ c-a\succ(b\succ c)=  (b\star a)\succ c-b\succ(a\succ c)\;.
\end{equation}%
Similarly, substituting $\a0,\,\b1$ and $\c0$ in (\ref{eq:L_dend_sat}), we get
\begin{equation}\label{eq:L-dend2}
  (a\succ b)\prec c-a\succ(b\prec c)=  (b\prec a)\prec c-b\prec(a\star c)\;.
\end{equation}
Since (\ref{eq:pre-Lie}) is invariant under the transposition which permutes $a$ and $b$, we have obtained all
relations, and so the relations for a $\cC_{pL}$-dendriform algebra are given by (\ref{eq:L-dend1}) and
(\ref{eq:L-dend2}). In the literature, such dendriform algebras are known as \emph{L-dendriform algebras} (see
\cite{Ldend}, where they have been introduced). One should keep in mind that, from our point of view, the L in
L-dendriform stands for \emph{pre-Lie}.
\end{example}
\begin{example}\label{exa:A3_dendri}
The defining relation for an $A_3$-associative algebra $(A,\mu)$ is 
\begin{equation}\label{eq:A_3_assoc_def}
  (ab)c+(bc)a+(ca)b=a(bc)+b(ca)+c(ab)\;.
\end{equation}%
It can be written in terms of associators in the following compact form:
\begin{equation}\label{eq:A_3_assoc_alt_def}
  \sum_{\s\in A_3} (a_{\s(1)},a_{\s(2)},a_{\s(3)})=0\;,
\end{equation}%
where $a_1,a_2,a_3\in A$. The symmetric form of (\ref{eq:A_3_assoc_alt_def}) is at the origin of the terminology
``$A_3$'' (see~\cite{gozeremm}); this form is often useful in computations, as we will see below.
Since (\ref{eq:A_3_assoc_alt_def}) is invariant under a cyclic permutation of $a_1,a_2,a_3$, the corresponding
dendriform algebras, which we will call \emph{$A_3$-dendriform algebras}, need to satisfy only one relation. We
obtain it by substituting $\a0,\,\b0$ and $\c1$ for $(a_1,x_1),\,(a_2,x_2)$ and $(a_3,x_3)$, in the relation
\begin{equation*}
 \sum_{\s\in A_3} ((a_{\s(1)},x_{\s(1)}),(a_{\s(2)},x_{\s(2)}),(a_{\s(3)},x_{\s(3)}))_{{}_\boxtimes}=0\;,
\end{equation*}%
where $(\cdot,\cdot,\cdot)_{{}_\boxtimes}$ stands for the associator of the product $\boxtimes$.  The resulting
relation defining $A_3$-dendriform algebras is given by
\begin{equation}\label{eq:A_3_first}
  a\succ(b\succ c)-(c\prec a)\prec b+c\prec(a\star b)=(a\star b)\succ c-b\succ (c\prec a)+(b\succ c)\prec a\;.
\end{equation}%
Notice that, upon defining $a\circ b:=a\succ b-b\prec a$ for all $a,b\in A$, the latter relation can be rewritten
in the following simple form:
\begin{equation}\label{eq:A_3_compact}
  (a\star b)\circ c-b\circ(c\prec a)-a\circ(b\succ c)=0\;.
\end{equation}%
We determine for this case also the relations of the corresponding tridendriform algebras. To do this, we need to
substitute in $\cR_\boxtimes=0$ at least one element from $A_1$ and the other ones from $A_0$. Notice that, if one
substitutes only one element from $A_1$, one obtains exactly the dendriform relations, with $\star$ standing now
for $a\star b:=a\prec b+a\succ b+a.b$, so these relations do not have to be computed again. Also, as above, there
is only one relation obtained by substituting two elements from $A_1$ and one from $A_0$, namely
\begin{equation}\label{eq:A_3_tridendri_2}
  (a.b)\prec c+(b\prec c).a+(c\succ a).b=a.(b\prec c)+b.(c\succ a)+c\succ(a.b)\;.
\end{equation}%
A final relation is obtained by substituting three elements from $A_1$. It is clear that the found relation just
says that $(A,.)$ is $A_3$-associative.
\end{example}
\begin{example}\label{exa:LA}
A \emph{Lie-admissible algebra} (or \emph{LA-algebra}) is classically defined as an algebra $(A,\mu)$ for which the
anticommutative product $\LB$, defined as the commutator $[a,b]:=ab-ba$, is a Lie bracket, i.e., satisfies the
Jacobi identity. The trilinear relation which characterizes Lie-admissible algebras is therefore given~by
\begin{equation}\label{equ:LA_def}
  \sum_{\s\in A_3}\left((a_{\s(1)},a_{\s(2)},a_{\s(3)})-(a_{\s(2)},a_{\s(1)},a_{\s(3)})\right)=0\;.
\end{equation}%
%
%
The relation (\ref{equ:LA_def}) is invariant under the full symmetry group $S_3$, so the corresponding dendriform
algebras, which we call \emph{LA-dendriform algebras}, are defined by a single relation, as in the case of
$A_3$-dendriform algebras. It is obtained in the same way as in that case, and is given by
\begin{eqnarray}
  \lefteqn{a\succ(b\succ c-c\prec b)-(b\succ c-c\prec b)\prec a-b\succ(a\succ c-c\prec a)}\label{LAcond}\\
  &+&(a\succ c-c\prec a)\prec b+c\prec(a\star b-b\star a)-(a\star b-b\star a)\succ c=0\;,\nonumber
\end{eqnarray}%
where $\star$ stands again for the sum of $\prec$ and $\succ$. As above, we define $a\circ b:=a\succ b-b\prec a$
for all $a,b\in A$ and observe that $a\star b-b\star a=a\circ b-b\circ a$, for all $a,b\in A$. Then the relation
defining $LA$-dendriform algebras can be rewritten in the following simple form:
\begin{equation}\label{eq:LA_compact}
  a\circ(b\circ c)-b\circ(a\circ c)-(a\circ b-b\circ a)\circ c=0\;.
\end{equation}%
It is equivalent to saying that $(A,\circ)$ is a pre-Lie algebra (see Example \ref{ex:L-dendri}).
\end{example}
\begin{example}\label{exa:AA}
An \emph{associative-admissible algebra} (or \emph{AA-algebra}) is similarly defined as an algebra $(A,\mu)$ for
which the commutative product $\LB^+$, defined as the anticommutator $\lb{a,b}^+:=ab+ba$, is associative. They are
in a certain sense the commutative analogs of LA-algebras. AA-algebras are characterized by the trilinear relation
\begin{equation}\label{equ:AA_def}
  (ab+ba)c+c(ab+ba)=a(bc+cb)+(bc+cb)a\;.
\end{equation}%
The relation (\ref{equ:AA_def}) is again invariant under the full symmetry group $S_3$, so the corresponding
dendriform algebras, \emph{AA-dendriform algebras}, are defined by a single relation. It is most easily obtained
from the compact form $\lb{\lb{a,b}^+,c}^+=\lb{a,\lb{b,c}^+}^+$ of the relation (\ref{equ:AA_def}). Indeed, let us
denote by $\LB^+_\boxtimes$ the anticommutator of $\boxtimes$, and let $a*b:=a\succ b+b\prec a$ for all $a,b\in A$
(not to be confused with $a\star b=a\succ b+a\prec b$). Using the obvious identity $a\star b+b\star a=a*b+b*a$ it
is easy to derive from Table \ref{tab:boxtimes} that
\begin{equation*}
  \lb{\a0,\b0}_\boxtimes^+=\ul{a*b+b*a}0\;,\qquad  \lb{\a0,\b1}_\boxtimes^+=\ul{a*b}1\;,
\end{equation*}%
for $a,b\in A$. Substituted in $\lb{\lb{\a0,\b0}^+_\boxtimes,\c1}^+_\boxtimes=
\lb{\a0,\lb{\b0,\c0}^+_\boxtimes=}^+_\boxtimes$, we obtain the following relation for AA-dendriform algebras:
\begin{equation}\label{eq:AA_compact}
  (a*b+b*a)*c=a*(b*c)\;.
\end{equation}%
This property is known as the \emph{Zinbiel} property, see \cite{lodaydend}.
\end{example}
\begin{example}\label{exa:P-algebra}
Our last example is closely related to Poisson algebras (see Examples \ref{exa:poisson} and
\ref{exa:P_to_Poisson}). Consider the following relation:
\begin{equation}\label{equ:poisson_pol}%
  3(ab)c=3a(bc)+(ac)b+(bc)a-(ba)c-(ca)b\;. 
\end{equation}%
In view of the mentioned relation to Poisson algebras, we call any algebra satisfying this relation a
\emph{P-algebra}. The category of all P-algebras is denoted by $\cP$. It was shown in \cite{GR2008} that P-algebras
are $A_3$-associative. Since (\ref{equ:poisson_pol}) admits no symmetry (when the variables $a,b,c$ are permuted),
we get three relations for the corresponding dendriform algebras, which we call {\em P-dendriform algebras}.  They
are given by the following formulas, where the first one is obtained by substituting $\a0,\,\b0$ and $\c1$ for
$a,\,b$ and~$c$ in (\ref{equ:poisson_pol}), where the product $\mu$ has been replaced by $\boxtimes$, and similarly
for the other two, where one substitutes $\a0,\,\b1,\,\c0$ and $\a1,\,\b0,\,\c0$ respectively:
\begin{align}
3\,(a\star b)\succ c&=3\,a\succ (b\succ c)+(a\succ c)\prec b+(b\succ c)\prec a \nonumber \\
&\quad-\;(b\star a)\succ c-(c\prec a)\prec b\;, \label{PD1}\\
3\,(a\succ b)\prec c&=3\,a\succ (b\prec c)+(a\star c)\succ b+(b\prec c)\prec a \nonumber \\
&\quad-\;(b\prec a)\prec c-(c\star a)\succ b\;, \label{PD3}\\
3(a\prec b)\prec c&=3\,a\prec (b\star c)+(a\prec c)\prec b+(b\star c)\succ a \nonumber \\
&\quad-\;(b\succ a)\prec c-(c\succ a)\prec b\;. \label{PD2}
\end{align}
In these formulas, $\star$ stands again for the sum of $\prec$ and $\succ$.
%
\end{example}
%
%

\subsection{Commutative and anticommutative dendriform algebras}\label{par:comm_anticom}
Many algebras of interest are commutative or anticommutative, i.e., they satisfy the relation $ab=ba$ or $ab=-ba$,
besides satisfying some other relations. It follows at once from the defining relations that:
\begin{enumerate}
  \item[(1)] Associative, pre-Lie, AA and P-algebras which are commutative, are precisely commutative
    associative algebras;
  \item[(2)] $A_3$-associative and  LA-algebras which are commutative, are just arbitrary (commutative)
    algebras; similarly, AA-algebras which are anticommutative are arbitrary (anticommutative) algebras:
  \item[(3)] $A_3$-associative, pre-Lie, LA and P-algebras which are anticommutative, are precisely Lie algebras;
  \item[(4)] Associative algebras which are anticommutative, are precisely (left and right) 2-step nilpotent
    algebras, i.e., algebras $A$ satisfying $(ab)c=a(bc)=0$ for all $a,b,c\in A$.
\end{enumerate}
It is clear from (\ref{eq:semi-direct}) that the corresponding dendriform algebras must satisfy the relation
$a\prec b=b\succ a$, respectively $a\prec b=-b\succ a$. It leads to the following definition.
\begin{defn}
  A $\cC$-dendriform algebra $(A,\prec,\succ)$ is said to be \emph{commutative} (resp.\ \emph{anticommutative}) if
  it satisfies $b\succ a=a\prec b$ (resp.\ $b\succ a=-a\prec b$) for all $a,b\in A$.
\end{defn}
In these cases it is natural to view $A$ as an algebra with only one product, by setting for all $a,b\in A$,
$a\times b:=a\succ b$ and the relations which $\times$ has to satisfy follow easily by substituting in the already
found dendriform relations everywhere $a\times b$ for $a\succ b$ and for $\pm b\prec a$, the sign depending on
whether commutativity or anticommutativity is considered. We give a few examples, based on the examples from
Section \ref{par:examples}.
\begin{example}\label{exa:a_com}
We start with (1) above: to obtain the relations of a commutative associative dendriform algebra, we substitute
$a\times b$ for $a\succ b$ and for $b\prec a$ in the relations (\ref{dend1}) -- (\ref{dend3}), to find the
relations
\begin{equation}\label{equ:dendri_com}
  (a\times b+b\times a)\times c=a\times(b\times c)\;, \quad c\times(a\times b)=a\times(c\times b)\;.
\end{equation}%
The first property is the {Zinbiel} property (see Example \ref{exa:AA}). The second property is know as the
\emph{NAP} (for non-associative, permutative) property, see \cite{livernet}.  Since the Zinbiel property implies
the NAP property, commutative associative dendriform algebras are, written in terms of a single product, the
same as Zinbiel algebras.
\end{example}
\begin{example}
  For (2) above, arbitrary (anti-) commutative algebras, one only gets the dendriform relation $a\prec b=\pm b\succ
  a$, with no relation for~$\times$.
\end{example}
\begin{example}
  For Lie algebras (case (3) above), the quickest way to obtain the relation which $\times$ must satisfy is by
  substituting $2a\times b$ (or just $a\times b$) for $a\circ b$ in (\ref{eq:LA_compact}), so we get the pre-Lie
  relation (\ref{eq:pre-Lie}). Thus, Lie dendriform algebras are, written in terms of a single product, pre-Lie
  algebras.
\end{example}
\begin{example}
  By definition, (right and left) $2$-step nilpotent algebras (case (4) above) satisfy $(ab)c=a(bc)=0$. Their
  dendriform algebras satisfy the following six relations:
  \begin{eqnarray*}
    (a\prec b)\succ c=(a\succ b)\prec c= (a\prec b)\prec c=0\;,\\
    c\prec(b\succ a) =c\succ(b\prec a)= c\succ(b\succ a)=0\;.
  \end{eqnarray*}
  It follows that anticommutative associative dendriform algebras are, in terms of a single product, also (right
  and left) $2$-step nilpotent algebras, as they satisfy the relation $(a\times b)\times c=a\times(b\times c)=0$.
\end{example}
\begin{remark}
Similarly, a tridendriform algebra is said to be \emph{commutative} or \emph{anticommutative} if it satisfies the
relations $a\succ b=\pm b\prec a$ and $a.b=\pm b.a$, with the plus sign of course corresponding to the commutative
case. Such tridendriform algebras are naturally seen as algebras with two operations ``$\times$'' and ``$.$'', upon
setting $a\times b:=a\succ b$, while keeping ``$.$''.
\end{remark}

\begin{example}
We give an example of an anticommutative tridendriform algebra: a Lie tridendriform algebra. We obtain the
relations from the relations of an $A_3$-tridendriform algebra, given in Example \ref{exa:A3_dendri}, by replacing
in them $a\succ b$ and $-b\prec a$ by $a\times b$, in particular $a\star b$ by $a\times b-b\times a+a.b$ and
$a\circ b$ by $2a\times b$, and using that $a.b=-b.a$. After some trivial simplifications, one finds that a Lie
tridendriform algebra is a Lie algebra, satisfying the following two relations, obtained from
(\ref{eq:A_3_compact}) and~(\ref{eq:A_3_tridendri_2}):
\begin{align*}
  (a.b)\times c&=a\times (b\times c)-(a\times b)\times c-b\times (a\times c)+(b\times a)\times c\;,\\
  c\times(a.b)&=(c\times b).a-(c\times a).b\;.
\end{align*}
In the literature, Lie tridendriform algebras are known as Post-Lie algebras (see \cite{baiguoni,vallette}).
\end{example}
\subsection{Categories of generalized dendriform algebras}\label{par:relating_examples}
Let, as before, $\cR_1=0,\dots,\cR_k=0$ be given relations. Recall that we denote by~$\cC$ the category of all
algebras $(A,\mu)$ over $\RF$ which satisfy these relations, with algebra homomorphisms as morphisms in
$\cC$. Clearly, the class of all $\cC$-dendriform algebras (over $\RF$) also form a category $\catdend$, with
morphisms the algebra homomorphisms. For example, the category of Loday dendriform algebras (constructed from
associative algebras) is denoted by $\hbox{Assoc}^{\hbox{\tiny dend}}$ and the category of $P-$dendriform algebras
is denoted by $\cP^{\hbox{\tiny dend}}$.

By the above, $\catdend$ is constructed out of $\cC$, but that does not mean that we know how to associate to
algebras in $\cC$ dendriform algebras in $\catdend$; we have on the contrary a (faithful) functor $\catdend\to\cC$,
which on objects $(A,\prec,\succ)$ is defined by $(A,\prec,\succ)\mapsto (A,\star)$, where $\star$ denotes, as in
the case of a Loday dendriform algebra, the sum of the products $\prec$ and~$\succ$; on morphisms, the functor is
just the identity in the sense that it sends the map underlying a morphism to itself.

Suppose that we have a second collection of relations $\cR'_1=0,\dots,\cR'_\ell=0$, where every $\cR_i$ is a linear
combination of $\cR'_1,\dots,\cR'_\ell$. It is clear that every algebra satisfying all relations
$\cR'_1=0,\dots,\cR'_\ell=0$ satisfies all relations $\cR_i=0$, and so $\catp$, the category of all algebras
satisfying the relations $\cR'_1=0,\dots,\cR'_\ell=0$, is a subcategory of $\cC$. Then $\catpdend$ is a subcategory
of $\catdend$, since the relations $\cR'_i=0$ can be seen as a subset of the relations $\cR_j=0$, and similarly for
the dendriform relations obtained from the relations $\cR'_i=0$ and $\cR_j=0$. Thus, we have the following
commutative diagram of categories:
\begin{center}
  \begin{tikzcd}[row sep=8ex, column sep=15ex]
    \cC',\mu'\arrow[hook]{r}{}&\cC,\mu\\ \catpdend,\prec',\succ'\arrow{u}{a\prec' b+a\succ' b}
    \arrow[hook]{r}{}&\catdend,\prec,\succ\arrow[swap]{u}{a\prec b+a\succ b}
  \end{tikzcd}
\end{center}
In this diagram, the horizontal arrows are inclusions and the products denote typical products of the objects of
the respective categories.

As a first application, we denote by $\catcom$ (resp.\ by $\catdendcom$) the subcategory of $\cC$ (resp.\ of
$\catdend$) consisting of all commutative algebras in the respective category. Then we have the following
commutative diagram of categories:
\begin{center}
  \begin{tikzcd}[row sep=8ex, column sep=15ex]
    \catcom,\mu\arrow[hook]{r}{}&\cC,\mu\\ \catdendcom,\prec,\succ\arrow{u}{a\prec b+a\succ b}
    \arrow[hook]{r}{}&\catdend,\prec,\succ\arrow[swap]{u}{a\prec b+a\succ b}
  \end{tikzcd}
\end{center}
Indeed, we can view the commutative algebras in $\cC$ as being those which satisfy the extra condition of
commutativity, and this relation leads to the condition of commutativity for the corresponding $\cC$-dendriform
algebras, by the above observation. The same applies, of course, to anticommutative algebras.
%
%

As a second application, we show how the above examples of $\cC$-dendriform algebras are related. We have the
following strict inclusion relations between the original category of algebras on the left; they lead to inclusion
relations between their corresponding categories of dendriform algebras on the right.
\begin{center}
  \begin{tikzcd}[row sep=9ex, column sep=3ex]
    \cP\arrow[hook]{d}&&\hbox{Assoc}\arrow[hook]{dll}\arrow[hook]{d}&\catPdend\arrow[hook]{d}
    &&\hbox{Assoc}^{\hbox{\tiny dend}}\arrow[hook]{dll}\arrow[hook]{d}\\
    \hbox{$A_3$-assoc}\arrow[hook]{dr}&&\hbox{pre-Lie}\arrow[hook]{dl}&
      \hbox{$A_3$}^{\hbox{\tiny dend}}\arrow[hook]{dr}&&\hbox{L}^{\hbox{\tiny dend}}\arrow[hook]{dl}\\
    &\hbox{LA}&&&\hbox{LA}^{\hbox{\tiny dend}}
  \end{tikzcd}
\end{center}
We have not included AA-algebras and their dendriform algebras, because there are no apparent inclusion relations
between the category of AA-algebras and any of the other categories that we considered.

The following table shows that the induced inclusions in the rightmost diagram are also strict and that there is no
inclusion relation between $\hbox{$A_3$}^{\hbox{\tiny dend}}$ or $\catPdend$ and $\hbox{L}^{\hbox{\tiny dend}}$. In
the table, the algebra $(A,\prec,\succ)$ is a free module of rank at least two and $a$ and $b$ are elements of a
basis of $A$. The first two columns describe the products $\prec$ and $\succ$ on some of the basis elements; it
is understood that all other products between elements of the basis are zero.

\medskip

\begin{table}[h]
  \def\arraystretch{1.6}
  \setlength\tabcolsep{0.2cm}
\centering
\begin{tabular}{|c||c||c||c||c|}
  \hline
   $\prec$&$\succ$&$\star$&\hbox{of type}&\hbox{not of type}\\
  \hline  \cline{1-5}
  $a\prec a=-b$&$\begin{array}{cc}a\succ a=a+b\\ b\succ a=b\end{array}$&
    $\begin{array}{cc} a\star a=a\\ b\star a=b\end{array}$&$A_3$-dendri&
    $\begin{array}{cc}\hbox{L-dendri}\\ \hbox{P-dendri}\end{array}$\\
  \hline 
  $a\prec b=b$&$\begin{array}{cc}b\succ a=b\\ b\succ b=b\end{array}$&
    $\begin{array}{cc}a\star b=b\\ b\star a=b\\ b\star b=b\end{array}$&
    \hbox{LA-dendri}&
    $\begin{array}{cc}\hbox{$A_3$-dendri}\\ \hbox{L-dendri}\end{array}$\\
  \hline
  ---&$b\succ a=a$&$b\star a=a$&L-dendri&
    $\begin{array}{cc}\hbox{$A_3$-dendri}\\ \hbox{dendri}\end{array}$\\
  \hline
  $a\prec b=-a$&$b\succ a=a$&$\begin{array}{cc} a\star b=-a\\ b\star a=a\end{array}$&P-dendri&L-dendri\\
  \hline
  $a\prec a=a+b$&---&$a\star a=a+b$&dendri&P-dendri\\
  \hline
\end{tabular}
\medskip
\caption{Some examples of generalized dendriform algebras.}\label{tab:examples}
\end{table}


%% file: rb.tex
\section{(Weak) Rota-Baxter operators}\label{sec:A3}
In this section, we introduce the notion of a weak Rota-Baxter operator, which generalizes the notion of a
Rota-Baxter operator. We show how such operators can be used to construct generalized dendriform algebras and give
an application to coboundary $\epsilon$-bialgebras.

\subsection{Dendriform algebras from Rota-Baxter operators}\label{par:dendri_from_RB}
We start with the definition of a Rota-Baxter operator (on an arbitrary algebra), see \cite{Guo}.
\begin{defn}
  Let $(A,\mu)$ be any algebra, let $\fR:A\to A$ be a linear map and let $\l\in R$. One says that $\fR$ is a
  \emph{Rota-Baxter operator of weight $\l$} of $A$ if $\fR$ satisfies the \emph{Rota-Baxter equation}
  \begin{equation}\label{eq:rota-baxter}
    \fR(a\fR(b)+\fR(a)b+\l ab)-\fR(a)\fR(b)=0\;,
  \end{equation}%
for all $a,b,\in A$. When $\l=0$ one simply speaks of a \emph{Rota-Baxter operator}. 
\end{defn}
Let $\cC$ be the category of all algebras satisfying a given collection of multilinear relations
$\cR_1=0,\dots,\cR_k=0$. We show in the following proposition how any Rota-Baxter operator (of weight zero) on any
algebra $(A,\mu)$ of $\cC$ leads to a $\cC$-dendriform algebra $(A,\prec,\succ)$.

\begin{prop}\label{prp:dendri_from_RB}
  Let $\fR$ be a Rota-Baxter operator on an algebra $(A,\mu)$ which belongs to $\cC$. For $a,b\in A$, let $a\succ
  b:=\fR(a)b$ and $a\prec b:=a\fR(b)$. Then $(A,\prec,\succ)$ is a $\cC$-dendriform algebra.
\end{prop}
\begin{proof}
We will give the proof for a trilinear relation $\cR=0$; it is easily generalized to $n$-linear relations by
induction on $n$. Recall from Section~\ref{sec:dendri} that $\cR=0$ leads to 3 dendriform relations which are
obtained by substituting two elements from~$A_0$ and one element from $A_1$ in $\cR_\boxtimes=0$, where $\boxtimes$
is the product on $A\times A$, defined by (\ref{eq:box_def}). Recall also that we write $\a0$ for $(a,0)$ and $\a1$
for $(0,a)$, where $a\in A$.

We show that such a substitution in $\cR_\boxtimes$ amounts to writing $\cR$ for three elements of $A$, on two of
which $\fR$ has been applied, and rewriting the result in terms of the dendriform operations. To show this, we
compare the effect of these substitutions on the two types of monomials $(ab)c$ and $a(bc)$, where each time we
consider the three possible substitutions. In view of Table \ref{tab:boxtimes}, the definition of $\prec$ and
$\succ$, and the Rota-Baxter equation (\ref{eq:rota-baxter}), we get for the first type the following correspondence:
\begin{equation*}
  \renewcommand{\arraystretch}{1.5}
  \begin{array}{rclrcl}
  (\a1\boxtimes\b0)\boxtimes\c0&=&\ul{(a\prec b)\prec c}1&=&\ul{(a\fR(b))\fR(c)}1\;,\\
  (\a0\boxtimes\b1)\boxtimes\c0&=&\ul{(a\succ b)\prec c}1&=&\ul{(\fR(a)b)\fR(c)}1\;,\\
  (\a0\boxtimes\b0)\boxtimes\c1&=&\ul{(a\star b)\succ c}1&=&\ul{(\fR(a)\fR(b))c}1\;,\\
  \end{array}
\end{equation*}
and similarly for the other type.  In the third line we have used  (\ref{eq:rota-baxter})
with $\l=0$, which says that $\fR:(A,\star)\to(A,\mu)$ is a morphism.
%
\end{proof}
\begin{remark}
Our proof shows that the $\cC$-dendriform relations can also formally be obtained from the relations
$\cR_i=0$ by formally applying $\fR$ to two of the variables and rewriting the resulting expression
in terms of the dendriform operations (using the Rota-Baxter equation). Our proof also explains where the
particular form of the Rota-Baxter equation comes from.
\end{remark}
%
%
As a direct consequence of Proposition \ref{prp:dendri_from_RB}, we have the following result, which is well-known
in the case of an associative or Lie algebra:
\begin{cor}\label{cor:RB_to_star}
  Let $\fR$ be a Rota-Baxter operator on an algebra $(A,\mu)$ in $\cC$. For $a,b\in A$, let $a\star
  b:=a\fR(b)+\fR(a)b$. Then $(A,\star)$ also belongs to $\cC$.
\end{cor}
\begin{remark}
The proof of Proposition \ref{prp:dendri_from_RB} is easily adapted to prove the following generalization of
Proposition \ref{prp:dendri_from_RB}: If $\fR$ is a Rota-Baxter operator of weight $\l$ on an algebra $(A,\mu)$
which belongs to $\cC$, then $(A,\prec,\succ,.)$ is a $\cC$-tridendriform algebra, upon defining $a\succ b:=\fR(a)b$
and $a\prec b:=a\fR(b)$ and $a.b:=\l ab$, for all $a,b\in A$. If fact, it suffices to change in the proof the
meaning of $a\star b$, which should now stand for $a\succ b+a\prec b+a.b$.
\end{remark}
\begin{remark}
In the case of Lie algebras, one encounters also the following equation, generalizing the Rota-Baxter equation (of
weight zero):
  \begin{equation}\label{eq:m-yang-baxter}
    \fR(a\fR(b)+\fR(a)b)=\fR(a)\fR(b)+\nu ab\;,
  \end{equation}%
where $\nu\in R$ is a constant. Equation (\ref{eq:m-yang-baxter}) is known as the \emph{modified Yang-Baxter
equation} and has many application in the theory of integrable systems (see \cite[Sect.\ 4.4.3]{AMV}). The
statement and proof of Proposition \ref{prp:dendri_from_RB}, and hence also Corollary \ref{cor:RB_to_star},
generalize easily to this case, so if $(A,\mu)\in\cC$ is equipped with a solution $\fR$ of the modified Yang-Baxter
equation (\ref{eq:m-yang-baxter}) then $(A,\prec,\succ)$ is a $\cC$-dendriform algebra, where $a\succ b:=\fR(a)b$
and $a\prec b:=a\fR(b)$ for all $a,b\in A$. It is clear that in the display in the proof of the proposition, we
only need to replace in line $3$, $\fR(a)\fR(b)$ by $\fR(a)\fR(b)+\nu ab$. For the rest
the proof is unchanged: these extra terms will disappear because the original product $\mu$ satisfies the relation
$\cR=0$.
\end{remark}

\begin{example}
The prime example of a solution to the modified Yang-Baxter equation is based on the notion of a Lie algebra
splitting (see \cite[Sect.\ 4.4.1]{AMV}). It naturally generalizes as follows. Let $\cC$ be, as before, the
category of all algebras satisfying a given set of relations. A \emph{$\cC$-algebra splitting} of $(A,\mu)\in\cC$
is a module direct sum decomposition $A=A_+\oplus A_-$ of $A$, where $A_+$ and $A_-$ are subalgebras of $A$. If one
denotes by $P_+$ and $P_-$ projection on $A_+$ and $A_-$, then $\fR:=P_+-P_-$ is a solution to
(\ref{eq:m-yang-baxter}), with $\nu=1$. Indeed, upon setting $a_+:=P_+(a)$ and $a_-:=P_-(a)$ for $a\in A$, one has,
for any $a,b\in A$,
\begin{equation*}
  \fR(a\fR(b)+\fR(a)b)=2\fR(a_+b_+-a_-b_-)=2(a_+b_++a_-b_-)\;,
\end{equation*}%
where we have used in the last step that $A_+$ and $A_-$ are subalgebras of $A$; this is clearly equal to
\begin{equation*}
  (a_+-a_-)(b_+-b_-)+(a_++a_-)(b_++b_-)=\fR(a)\fR(b)+ab\;.
\end{equation*}%
It follows that a $\cC$-algebra splitting of $(A,\mu)\in\cC$ yields an algebra $(A,\star)\in\cC$, where $a\star
b:=a\fR(b)+\fR(a)b=a(b_+-b_-)+(a_+-a_-)b$, for $a,b\in A$.
\end{example}
\subsection{Dendriform algebras from weak Rota-Baxter operators}\label{par:dendri_from_weak_RB}
We now introduce the notion of a \emph{weak} Rota-Baxter operator, which generalizes the notion of a Rota-Baxter
operator. For any algebra $(A,\mu)$, we denote by $C(A)$ the set of elements $c$ of $A$ which commute with all
elements in $A$. It is a submodule of $A$ but is in general not a subalgebra\footnote{For a general algebra, $C(A)$
strictly contains the center $Z(A)$, whose elements are required to have the extra property that any associator
containing them vanishes.}
of $A$.

\begin{defn}
  Let $\fR:A\to A$ be a linear map and let $\l\in R$. One says that $\fR$ is a \emph{weak Rota-Baxter operator of
  weight $\l$} of $A$ if, for all $a,b\in A$, 
  \begin{equation}\label{eq:weak_rb}
    \fR(a\fR(b)+\fR(a)b+\l ab)-\fR(a)\fR(b)\in C(A)\;.
  \end{equation}%
  When $\l=0$ one simply speaks of a \emph{weak Rota-Baxter operator} of $A$. 
\end{defn}
We show how Proposition \ref{prp:dendri_from_RB} can be generalized to the case of weak Rota-Baxter operators (of
weight zero). For clarity, and in view of the examples, we will restrict ourselves to the case of trilinear
relations. We say that a trilinear relation $\cR=0$ has \emph{commutator form} if it can be written as a linear
combination of terms which have the form $\lb{ab,c}=(ab)c-c(ab)$ $(=-\lb{c,ab})$. Said differently, $\cR=0$ has
{commutator form} if $\cR$ is of the form
\begin{equation}
  \cR(a_1,a_2,a_3)=\sum_{\s\in S_3}c_\s\lb{a_{\s(1)}a_{\s(2)},a_{\s(3)}}\;,
\end{equation}%
for some constants $c_\s\in R$. A set of trilinear relations $\cR_1=0,\dots,\cR_k=0$, is said to have
\emph{commutator form} if it can equivalently be written as a set of trilinear relations, where each
relation has commutator form.

\begin{prop}\label{prp:dendri_from_weak_RB}
  Let $\cR_1=0,\dots,\cR_k=0$ be a collection of trilinear relations which are assumed to have commutator form. Let
  $\cC$ be the category of all algebras satisfying these relations. Let $\fR$ be a weak Rota-Baxter operator on an
  algebra $(A,\mu)$ which belongs to $\cC$. For $a,b\in A$, define $a\succ b:=\fR(a)b$ and $a\prec
  b:=a\fR(b)$. Then $(A,\prec,\succ)$ is a $\cC$-dendriform algebra.
\end{prop}

\begin{proof}
By the assumption, we may assume that $\cR_1=0,\dots,\cR_k=0$ have commutator form. Let $\cR=0$ be one of these
relations. We can repeat for $\cR$ the proof of Proposition~\ref{prp:dendri_from_RB}, except that we need to show
how to express the terms of the form $(\fR(a)\fR(b))c$ and $c(\fR(a)\fR(b))$ in terms of the dendriform operations
and that by this procedure the same terms are obtained as by substituting in $\cR_\boxtimes$ two terms from $A_0$
and one term from $A_1$. To do this, first observe that (\ref{eq:weak_rb}) can (for $\l=0$) be equivalently written
as the condition that $\lb{\fR(a\star b),c}=\lb{\fR(a)\fR(b),c},$ where $a\star b=a\succ b+a\prec
b=a\fR(b)+\fR(a)b$, leading to the following correspondence:
\begin{equation*}
  \lb{\a0\boxtimes\b0,\c1}_\boxtimes=\ul{(a\star b)\succ c-c\prec(a\star b)}1=\ul{\lb{\fR(a)\fR(b),c}}1\;,
\end{equation*}
where $\LB_\boxtimes$ stands for the commutator of the product $\boxtimes$. For the two other possible
substitutions, it is not necessary to use the commutator form and one can simply rely on the formulas given in the
proof of Proposition~\ref{prp:dendri_from_RB}. Yet, for completeness, we also express them in commutator form:
\begin{eqnarray*}
  &&\lb{\a0\boxtimes\b1,\c0}_\boxtimes=\ul{(a\succ b)\prec c-c\succ(a\succ b)}1=\ul{\lb{\fR(a)b,\fR(c)}}1\;,\\
  &&\lb{\a1\boxtimes\b0,\c0}_\boxtimes=\ul{(a\prec b)\prec c-c\succ(a\prec b)}1=\ul{\lb{a\fR(b),\fR(c)}}1\;.
\end{eqnarray*}
It follows that the $\cC$-dendriform relation $\cR=0$ is satisfied by the products $\prec$ and $\succ$, defined by
the weak Rota-Baxter operator $\fR$.
\end{proof}
The above theorem can be applied to $A_3$-associative algebras and Lie admissible algebras,
since (\ref{eq:A_3_assoc_def}) and (\ref{equ:LA_def}) can be respectively rewritten in the commutator forms
\begin{align}
  \lb{ab,c}+\lb{bc,a}+\lb{ca,b}=0\;,\label{eq:A_3_assoc_comm_form}\\
  \sum_{\s\in A_3}\lb{a_{\s(1)} a_{\s(2)}-a_{\s(2)} a_{\s(1)},a_{\s(3)}}=0\;.\label{equ:LA_def_comm_form}
\end{align}
However, many relations cannot be written in commutator form. The associativity relation, $a(bc)=(ab)c$, is a prime
example; other examples are the derivation property $a(bc)=(ab)c+b(ac)$, the Zinbiel property $a(bc)=(ab+ba) c$ and
the NAP property $a(bc)=b(ac)$, just to mention a few. In such cases, when the relations of $\cC$ imply a relation
$\cR=0$ which can be written in commutator form, any dendriform algebra $(A,\prec,\succ)$ obtained by using a weak
Rota-Baxter operator on an algebra $(A,\mu)$ in $\cC$ will satisfy (at least) the $\cC$-dendriform relation, derived
from $\cR=0$. Moreover, any relation $\cR=0$ which does not involve a product of two of the variables leads to a
(single) dendriform relation. We illustrate this in the following example, on which we will elaborate in the
following subsection.

\begin{example}\label{exa:a3_weak}
The associativity relation, $a(bc)=(ab)c$ can clearly not be written in commutator form. Summing up three instances
of this relation it implies however $(ab)c+(bc)a+(ca)b=a(bc)+b(ca)+c(ab)$, which is the relation of
$A_3$-associativity, which we wrote in commutator form in (\ref{eq:A_3_assoc_comm_form}). Therefore, if $\fR$ is a
weak Rota-Baxter operator on an associative algebra $(A,\mu)$ then $(A,\prec,\succ)$, with $\prec$ and $\succ$
defined by $a\prec b:=a\fR(b)$ and $a\succ b:=\fR(a)b$ is a priori not a Loday dendriform algebra, but it is an
$A_3$-dendriform algebra. Moreover, the associativity relation $a(bc)=(ab)c$ does not contain a product of $a$ and
$c$, so we do not need to use the weak Rota-Baxter equation to rewrite $\fR(a)(b\fR(c))=(\fR(a)b)\fR(c)$ in terms
of the dendriform products. The resulting relation $a\succ(b\prec c)=(a\succ b)\prec c$ of $(A,\prec,\succ)$ is
called \emph{inner-associativity}.
\end{example}
It follows that a weak Rota-Baxter operator on an associative algebra leads to an inner-associative
$A_3$-dendriform algebra. We show in the following example that in general the latter is not a Loday dendriform
algebra.
\begin{example}
  Let $A$ be a commutative associative algebra. Every linear map $\cR:A\to A$ is a weak Rota-Baxter operator since
  $C(A)=A$, hence leads to an inner-associative $A_3$-dendriform algebra. To see that it may not be a classical
  dendriform algebra, take $\cR=\Id_A$. Then $a\prec b=a\succ b=ab$ and (\ref{dend1}) cannot be satisfied, unless
  $abc=0$ for all $a,b,c\in A$.
\end{example}
\begin{remark}
The proof of Proposition \ref{prp:dendri_from_weak_RB} is easily adapted to prove the following generalization of
Proposition \ref{prp:dendri_from_weak_RB}: under the same assumptions on the relations of $\cC$, any weak
Rota-Baxter operator $\fR$ of weight $\l$ on an algebra $(A,\mu)\in\cC$ leads to a
$\cC$-tridendriform algebra, upon setting $a\succ b:=\fR(a)b$ and $a\prec b:=a\fR(b)$ and $a.b:=\l ab$, for all
$a,b\in A$. Again, it suffices to change in the proof the meaning of $a\star b$, which should now stand for $a\succ
b+a\prec b+a.b$. The comments made about relations which cannot be written in commutator form apply here without
modification.
\end{remark}
\begin{remark}
If we denote by $C'(A)$ the set of elements $c$ of $A$ which anticommute with all elements of $A$, i.e., $ac=-ca$
for all $a\in A$, we can also consider operators $\cR$ satisfying (\ref{eq:weak_rb}), with $C(A)$ replaced by
$C'(A)$. The results of this section are easily adapted to the case of such operators. For example, the conclusion
of Proposition \ref{prp:dendri_from_weak_RB} still holds for such an operator $\cR$ when the relations have
anticommutator form. An example of such a relation is the relation (\ref{equ:AA_def})
defining AA-algebras.
\end{remark}

\subsection{Application: coboundary $\epsilon$-bialgebras}\label{par:coboundary}
As an application of weak Rota-Baxter operators, we now generalize  a result obtained by M. Aguiar in
\cite{aguiar_inf}, which we will recall.
We first recall the definition of an $\epsilon$-bialgebra:
\begin{defn}
An \emph{$\epsilon$-bialgebra} is a triple $(A,\mu,\Delta)$, where $A$ is  an $R$-mo\-dule and 
$\mu:A\otimes A\rightarrow A$ and $\Delta:A\to A\otimes A$ are linear maps, such that 
\begin{enumerate}
  \item[(1)] $\mu$ is associative;
  \item[(2)] $\Delta$ is coassociative;
  \item[(3)] $\Delta$ is a derivation: $\Delta(ab)=a\cdot\Delta(b)+\Delta(a)\cdot b$, for all $a,b\in A$.
\end{enumerate}
%
\end{defn}
In item (3), we have used a dot to denote the natural left, resp.\ right action of $A$ on $A\otimes A$; later on in
this section, it will also be used for the natural left and right actions of $A$ on $A\otimes A\otimes A$.

Let $(A,\mu,\Delta)$ be an $\epsilon$-bialgebra and let us write $\Delta(a)=\sum_{(a)}a_{(1)}\otimes a_{(2)}$ for
all $a\in A$ (Sweedler's notation). It is shown in \cite{aguiar_inf} that if one defines $\po ab:=\sum_{(b)}b_{(1)}ab_{(2)}$ for all $a,b\in A$, then $(A,\circ)$ is a pre-Lie algebra.
%
This yields a functor which associates to any $\epsilon$-bialgebra $(A,\mu,\Delta)$ the corresponding pre-Lie
algebra $(A,\circ)$, and which is identity on morphisms.

The fundamental observation of Aguiar is that the restriction of this functor to quasi-triangular
$\epsilon$-bialgebras factors in a natural way through the category of Loday dendriform algebras, as in the
following diagram:
\begin{center}
  \begin{tikzcd}[row sep=8ex, column sep=15ex]
    \hbox{QT\ }\epsilon\hbox{-bialg},\mu,r\arrow{d}{\sum_iau_ibv_i,\sum_iu_iav_ib}\arrow{r}{r\cdot a-a\cdot r}
      &\epsilon\hbox{-bialg},\mu,\Delta\arrow{d}{\sum_{(b)}b_{(1)}ab_{(2)}}\\
    \hbox{Assoc}^{\hbox{\tiny dend}},\prec,\succ\arrow{r}{a\succ b-b\prec a}&\hbox{pre-Lie},\circ
  \end{tikzcd}
\end{center}
In order to explain this diagram, we first recall from \cite{aguiar_inf} that a \emph{quasi-triangular
 $\epsilon$-bialgebra} is a triple $(A,\mu,r)$, where $(A, \mu )$ is an associative algebra and $r\in A\otimes A$
is a solution of the \emph{associative Yang-Baxter equation}
\begin{equation*}
  \AYB(r):=r_{13}r_{12}-r_{12}r_{23}+r_{23}r_{13}=0\;.
\end{equation*}%
Let $(A,\mu,r)$ be a quasi-triangular $\epsilon$-bialgebra and write $r$ as $r=\sum_i u_i\otimes v_i$. On the one
hand, setting for all $a\in A$
\begin{equation}\label{eq:delta_r_def}
  \Delta_r(a):=r\cdot a-a\cdot r\;, 
\end{equation}
we get an $\epsilon$-bialgebra $(A,\mu,\Delta_r)$. On the other hand, the map $\fR:A\rightarrow
A$, defined for all $a\in A$ by $\fR(a)=\sum _iu_iav_i$, is a Rota-Baxter operator for $A$, and so, by Proposition
\ref{prp:dendri_from_RB}, the products $\prec$ and $\succ$ defined for all $a,b\in A$ by
\begin{equation}\label{eq:A3_from_bia}
  a\prec b:=\sum_i au_ibv_i\;,\qquad a\succ b:=\sum_i u_iav_ib\;,
\end{equation}
make $(A,\prec,\succ$) into a Loday dendriform algebra.
%
The above construction which associates to a solution of the associative Yang-Baxter equation an
$\epsilon$-bialgebra has a natural generalization, given by the next proposition. It gives necessary and sufficient
conditions on $r\in A\otimes A$ so that the triplet $(A,\mu,\Delta_r)$ is an $\epsilon$-bialgebra, with $\Delta_r$
defined by (\ref{eq:delta_r_def}).
\begin{prop}[\cite{AguiarContemp}]\label{prp:cob_to_eps}
  Let $(A,\mu)$ be an associative algebra and let $r\in A\otimes A$. Then $(A,\mu,\Delta_r)$ is an
  $\epsilon$-bialgebra if and only if $\AYB(r)$ is \emph{invariant}, i.e., $a\cdot\AYB(r)=\AYB(r)\cdot a$, for all
  $a\in A$. One then says that $(A,\mu,r)$ is a \emph{coboundary $\epsilon$-bialgebra}.
\end{prop}
%


\begin{prop}\label{prp:cob_to_dendri}
  Let $(A,\mu,\sum_iu_i\otimes v_i)$ be a coboundary $\epsilon$-bialgebra. \\
(1) The linear map $\fR:A\to A$, defined for all $a\in A$ by $\fR(a):=\sum_i u_iav_i$, is a weak
      Rota-Baxter operator for $A$.\\
(2) For $a,b\in A$, let $a\succ b:=\fR(a)b=\sum_i u_iav_ib$ and $a\prec b:=a\fR(b)=\sum_i
      au_ibv_i$. Then $(A,\prec,\succ)$ is an inner-associative $A_3$-dendriform algebra.
\end{prop}
\begin{proof}
  We only need to prove (1), because (2) follows from it by Example~\ref{exa:a3_weak}.
  To do this, we show that the linear map $\omega:A\otimes A\to A$, defined for $a,b\in A$~by $\omega(a\otimes
  b):=\fR(a)\fR(b)-\fR(a\fR(b)+\fR(a)b)$ satisfies $\omega(a\otimes b)c=c\omega(a\otimes b)$ for all $a,b,c\in
  A$. We do this by relating $\omega$ with $\AYB(r)$. Without loss of generality, we may assume that the
  associative algebra $A$ has a unit, denoted~$1_A$. Writing $r=\sum_i u_i\otimes v_i$,
\begin{align*}
  \AYB(r)&=r_{13}r_{12}-r_{12}r_{23}+r_{23}r_{13}\\
         &=\sum_{i,j}(u_i\otimes1_A\otimes v_i)(u_j\otimes v_j\otimes1_A)-
            \sum_{i,j}(u_i\otimes v_i\otimes 1_A)(1_A\otimes u_j\otimes v_j)\\
         &\quad+\sum_{i,j}(1_A\otimes u_i\otimes v_i)(u_j\otimes 1_A\otimes v_j)\\
         &=\sum_{i,j}(u_iu_j\otimes v_j\otimes v_i-u_i\otimes v_iu_j\otimes v_j+u_j\otimes u_i\otimes v_iv_j)\;;\\
  \omega(a\otimes b)&=\sum_{i,j}u_iav_iu_jbv_j-\fR\(\sum_i au_ibv_i+\sum_iu_iav_ib\)\\
    &=-\sum_{i,j}(u_jau_ibv_iv_j+u_ju_iav_ibv_j-u_iav_iu_jbv_j)\\
    &=-\sum_{i,j}(u_iu_jav_jbv_i-u_iav_iu_jbv_j+u_jau_ibv_iv_j)\;.
\end{align*}
If we compare these two expressions and we write $\AYB(r)$ as $\AYB(r)=\sum_k X_k\otimes Y_k\otimes Z_k$,
then we see that $\omega(a\otimes b)=-\sum_k X_ka Y_kb Z_k$. The invariance of $\AYB(r)$, which can be written
as $\sum_k cX_k\otimes Y_k\otimes Z_k=\sum_k X_k\otimes Y_k\otimes Z_kc$ for all $c\in A$ therefore yields
$\omega(a\otimes b)c=-\sum_k X_ka Y_kb Z_kc=-\sum_k cX_ka Y_kb Z_k=c \omega(a\otimes b)$, as was to be shown.
\end{proof}
Proposition \ref{prp:cob_to_dendri} leads to the following commutative diagram, generalizing Aguiar's commutative
diagram:
\begin{center}
  \begin{tikzcd}[row sep=8ex, column sep=15ex]
    \hbox{Cob-}\epsilon\hbox{-bialg},\mu,r\arrow{d}{\sum_iau_ibv_i,\sum_iu_iav_ib}\arrow{r}{r\cdot a-a\cdot r}
      &\epsilon\hbox{-bialg},\mu,\Delta\arrow{d}{{\sum_{(b)}b_{(1)}ab_{(2)}}}\\
    \hbox{$A_3^{\hbox{\tiny dend}}$},\prec,\succ\arrow{r}{a\succ b-b\prec a}&\hbox{pre-Lie},\circ
  \end{tikzcd}
\end{center}


%

\subsection{Curved Rota-Baxter systems}
We show in this paragraph that curved Rota-Baxter systems also provide examples of inner-associative
$A_3$-dendri\-form algebras. We first recall the definition of such systems (see \cite{brzez}).
\begin{defn}
  Let $A$ be an associative algebra endowed with linear maps $\fR, \fS:A\rightarrow A$ and $\omega :A\otimes
  A\rightarrow A$. The 4-tuple $(A, \fR, \fS, \omega )$ is called a \emph{curved Rota-Baxter system} if the following
  conditions are satisfied, for all $a, b\in A$:
  \begin{align}
    \fR(a)\fR(b)&=\fR(\fR(a)b+a\fS(b))+\omega (a\otimes b)\;, \label{curved1} \\
    \fS(a)\fS(b)&=\fS(\fR(a)b+a\fS(b))+\omega (a\otimes b)\;. \label{curved2}
  \end{align}
\end{defn}
The definition is easily generalized to arbitrary algebras, but not the results which follow; this is why we
consider only the case of associative algebras. Notice that weak Rota-Baxter operators on an associative
algebra~$A$ correspond to curved Rota-Baxter systems $(A, \fR, \fS, \omega )$ with $\fR=\fS$ and having the
property that $\omega$ takes values in $Z(A)$, the center of $A$ (which coincides with $C(A)$ because $A$ is
associative). Under this correspondence, the following proposition generalizes item (2) of Proposition
\ref{prp:cob_to_dendri}.
\begin{prop}\label{curvedRB}
  Let $(A,\fR,\fS,\omega )$ be a curved Rota-Baxter system. Define two new products on $A$ by setting $a\succ
  b:=\fR(a)b$ and $a\prec b=a\fS(b)$, for all $a, b\in A$.  Then $(A,\prec,\succ )$ is an $A_3$-dendriform algebra
  if and only if $\omega$ takes values in $Z(A)$. In any case, $(A,\prec,\succ )$ is inner-associative.
\end{prop}
\begin{proof}
$(A,\prec,\succ)$ is inner-associative, since for all $a,b,c\in A$, 
\begin{equation*}
  (a\succ b)\prec c=(\fR(a)b)\prec c=\fR(a)b\fS(c)=a\succ (b\fS(c))=a\succ (b\prec c)\;. 
\end{equation*}
Using (\ref{curved2}) we find that
\begin{align*}
  (a\prec b)\prec c-a\prec (b\prec c+b\succ c)&=a\fS(b)\fS(c)-a\fS(b\fS(c)+\fR(b)c)\nonumber\\
  &=a\omega (b\otimes c)\;,
\end{align*}
and similarly, using (\ref{curved1}), $b\succ (c\succ a)-(b\prec c+b\succ c)\succ a=\omega (b\otimes c)a$.
So, (\ref{eq:A_3_first}) is satisfied (i.e., $(A,\prec,\succ)$ is an $A_3$-dendriform algebra) if and only if
$a\omega (b\otimes c)=\omega (b\otimes c)a$, for all $a, b, c\in A$; in turn, this is equivalent to $\omega
(b\otimes c)\in Z(A)$, for all $b,c\in A$.
\end{proof}

The proof also shows that when $\omega=0$ the $A_3$-dendriform algebra which is obtained is a Loday dendriform
algebra; this was already observed in \cite{Brzfirst}.

It was proven in \cite{brzez} that, if $(A,\fR,\fS,\omega)$ is a curved Rota-Baxter system and we define a new
product on $A$ by $a\circ b=\fR(a)b-b\fS(a)$, then $(A, \circ )$ is a pre-Lie algebra if and only if $\omega
(a\otimes b-b\otimes a)\in Z(A)$, for all $a, b\in A$. In particular, $(A,\circ)$ is a pre-Lie algebra when
$\omega$ takes values in $Z(A)$. We recover this result as a direct consequence of Example \ref{exa:A3_dendri} and
Proposition~\ref{curvedRB}.

%
%
\begin{example}\label{exa:rbs_gen}
Let $A$ be an associative algebra and $\fR, \fS:A\rightarrow A$ be a left, respectively right Baxter operator,
i.e., $\fR(a)\fR(b)=\fR(\fR(a)b)$ and $\fS(a)\fS(b)=\fS(a\fS(b))$, for all $a, b\in A$, satisfying the extra
condition that
\begin{equation*}
  \fR(a)\fS(b)=\fR(a\fS(b))=\fS(\fR(a)b)
\end{equation*}
for all $a,b\in A$. Then $(A,\fR,\fS,\omega)$ is a curved Rota-Baxter system, where $\omega :A\otimes A\rightarrow
A$ is defined by $\omega (a\otimes b)=-\fR(a)\fS(b)$. If moreover $\fR(a), \fS(a)\in Z(A)$ for all $a\in A$, then
$\omega$ takes values in $Z(A)$, hence Proposition \ref{curvedRB} can be applied to yield an (inner-associative)
$A_3$-dendriform algebra.  A particular case of this example already appears in \cite{brzez}, where it is shown
that if $r=\sum _ix_i\otimes y_i$ and $s=\sum_jz_j\otimes w_j$ are invariant, then the linear maps $\fR,
\fS:A\rightarrow A$ and $\omega :A\otimes A\rightarrow A$, defined for $a\in A$ by
\begin{equation*}
  \fR(a):=\sum _ix_iay_i\;,\qquad \fS(a):=\sum _jz_jaw_j\;,\qquad\omega (a\otimes b)=-\fR(a)\fS(b)\;, 
\end{equation*}
make $(A,\fR, \fS,\omega)$ into a curved Rota-Baxter system.
\end{example}


%% file: pol.tex
\section{Dendriform algebras in polarized form}\label{par:gendend2}

In this section, we introduce the notion of a dendriform algebra for algebras $(A,\cdot,\LB)$, where ``$\cdot$'' is
commutative and $\LB$ is anticommutative, satisfying again any finite collection of (extra) relations. It will be
shown that this notion of a dendriform algebra corresponds to the one introduced in Section~\ref{sec:dendri}, via a
polarization functor which we will also introduce. 

\subsection{Polarized algebras}
We first define the class of algebras which we will consider in this section.
\begin{defn}
  An algebra $(A,\cdot,\LB)$ is said to be a \emph{polarized algebra} when ``$\cdot$'' is commutative and $\LB$ is
  anticommutative, i.e., for all $a,b\in A$,
  \begin{equation*}
    b\cdot a=a\cdot b\;,\qquad\hbox{and}\qquad\lb{b,a}=-\lb{a,b}\;.
  \end{equation*}%
\end{defn}
The choice of the terminology \emph{polarized} will become clear in Section \ref{par:polarization}, when we will see
how we can obtain polarized algebras from algebras with one product by using a procedure called
\emph{polarization}. 

\begin{example}
If $(A,\cdot)$ is a commutative algebra, we can make it into a polarized algebra $(A,\cdot,\LB)$ simply by adding any
anticommutative product $\LB$ on $A$, for example the trivial (zero) product. Similarly, any anticommutative
algebra $(A,\LB)$ can be made into a polarized algebra.
\end{example}
\begin{example}\label{exa:poisson}
Recall (for example from \cite{PLV}) that an algebra $(A,\cdot,\PB)$ is a \emph{Poisson algebra} if $(A,\cdot)$ is
a commutative associative algebra, $(A, \PB)$ is a Lie algebra and the two products are compatible in the sense
that
\begin{equation*}
 \{a\cdot b, c\}=a\cdot\{b, c\}+\{a, c\}\cdot b\;,
\end{equation*}
for all $a, b, c\in A$. The latter condition can also be formulated by saying that the product $\PB$, usually
referred to as the \emph{Poisson bracket}, is a derivation in each one of its arguments.  Clearly, every Poisson
algebra $(A,\cdot,\PB)$ is a polarized algebra. We will come back several times to this example.
\end{example}

\subsection{Polarized $\cC$-dendriform algebras}\label{par:pol_dendri}
In analogy with Definition~\ref{def:dendri_gen}, we now define the notion of a dendriform algebra for a polarized
algebra.  Here, $\cR_1=0,\dots,\cR_k=0$ are given relations involving the products ``$\cdot$'' and $\LB$
(only). The category of all polarized algebras satisfying these relations is denoted by $\catpol$. The morphisms in
$\catpol$ are the algebra homomorphisms.
\begin{defn}
  An algebra $(A,*,\circ)$ is said to be a \emph{polarized $\cC$-dendriform algebra} if $\(A\times A,\odot,
  \GB\)\in\catpol$, where $\odot$ and $\GB$ are defined, for $(a,x)$ and $(b,y)$ in $A\times A$, by
  \begin{align}
    (a,x)\odot (b,y):=(a*b+b*a,a*y+b*x)\;,\label{equ:dendri_pol_1}\\
    \gb{(a,x), (b,y)}:=(a\circ b-b\circ a,a\circ y-b\circ x)\;.\label{equ:dendri_pol_2}
  \end{align}%
%
\end{defn}
The category of all polarized $\cC$-dendriform algebras (over $\RF$) is denoted by $\catdendpol$. The morphisms in
this category are the algebra homomorphisms. Setting $x=y=0$ in (\ref{equ:dendri_pol_1}) and in
(\ref{equ:dendri_pol_2}), we see that we have again a faithful functor $\catdendpol\to\catpol$, defined on objects
by $(A,*,\circ)\mapsto (A,\cdot,\LB)$, where the two new products on $A$ are defined, for all $a,b\in A$, by
\begin{equation}\label{eq:vert_pol}
  a\cdot b:=a*b+b*a\;,\qquad\hbox{and}\qquad \lb{a,b}:=a\circ b-b\circ a\;.
\end{equation}%
\begin{remark}
The above definition of a polarized $\cC$-dendriform algebra admits the following natural generalization.  An
algebra $(A,*,\circ,\mid\,,\square)$ is said to be a \emph{polarized $\cC$-tridendriform algebra} if
$(A,\mid\,,\square)$ is a polarized algebra and $\(A\times A,\odot, \GB\)\in\catpol$, where $\odot$ and $\GB$ are
defined for $(a,x)$ and $(b,y)$ in $A\times A$, by
\begin{align}
  (a,x)\odot (b,y)&:=(a*b+b*a+a\mid b,a*y+b*x+x\mid y)\;,\label{equ:tridendri_pol_1}\\
  \gb{(a,x), (b,y)}&:=(a\circ b-b\circ a+a\squ b,a\circ y-b\circ x+x\squ y)\;.\label{equ:tridendri_pol_2}
\end{align}%
We have a functor from the category $\cattridendpol$ of all polarized $\cC$-tridendriform algebras to $\catpol$,
defined on objects by $(A,*,\circ,\mid\,,\square)\mapsto (A,\cdot,\LB)$, where 
\begin{equation*}
  a\cdot b:=a*b+b*a+a\mid b\;,\qquad\hbox{and}\qquad \lb{a,b}:=a\circ b-b\circ a+a\squ b\;,
\end{equation*}%
for all $a,b\in A$. Any polarized $\cC$-dendriform algebra $(A,*,\circ)$ can be seen in a natural way as a
polarized $\cC$-tridendriform algebra by considering $(A,*,\circ,\vert,\squ)$, where the products $\vert$ and
$\squ$ are trivial.
\end{remark}
\subsection{Algebras defined by multilinear relations}\label{sec:multilinear_pol}
As in the case of a $\cC$-dendriform algebra (see Section \ref{sec:multilinear}), the relations which every
polarized $\cC$-dendriform algebra must satisfy, can be easily computed when the relations~$\cR_i=0$ of $\catpol$
are multilinear. We show again how these relations can be computed for a trilinear relation $\cR=0$. Thanks to
commutativity and anticommutativity, $\cR$ is of the form
\begin{align*}
  \cR(a_1, a_2, a_3)&=\sum_{\s\in A_3}\l_\s(a_{\s(1)}\cdot a_{\s(2)})\cdot a_{\s(3)}+\sum_{\s\in A_3}
      \l'_\s \lb{a_{\s(1)},\lb{a_{\s(2)},a_{\s(3)}}}\\
    &\quad+\sum_{\s\in A_3}\l''_\s\lb{a_{\s(1)},a_{\s(2)}}\cdot a_{\s(3)}+\sum_{\s\in A_3}\l'''_\s
      \lb{a_{\s(1)}\cdot a_{\s(2)},a_{\s(3)}}\;,
\end{align*}%
where the constants $\l_\s,\dots,\l_\s'''$ belong to the base ring $\RF$. Notice that we have the same number of
constants as in the case of an algebra with one product, namely~12; we will see the reason for this in Section
\ref{par:polarization}.

By trilinearity, the relations which must be satisfied by every algebra in $\catdendpol$ are obtained by demanding
that the relations are satisfied on all possible triplets of elements of $A\times A$, taken from the union of $A_0$
and $A_1$, which is a generating set of $A\times A$. In the following two tables, we exhibit the products $\odot$
and $\GB$ in terms of these generators:

\begin{table}[h]
  \def\arraystretch{1.8}
  \setlength\tabcolsep{0.4cm}
  \centering
\begin{tabular}{c|cc}
   $\odot$&$\b0$&$\b1$\\
  \hline 
  $\a0$&$\ul{a*b+b*a}0$&$\ul{a* b}1$\\
  $\a1$&$\ul{b*a}1$&$(0,0)$\\
\end{tabular}\quad
\begin{tabular}{c|cc}
   $\GB$&$\b0$&$\b1$\\
  \hline 
  $\a0$&$\ul{a\circ b-b\circ a}0$&$\ul{a\circ b}1$\\
  $\a1$&$-\;\ul{b\circ a}1$&$(0,0)$\\
\end{tabular}
\medskip
\caption{The products $\odot$ and $\GB$ for generators of $A\times A$.}\label{tab:odot_gb}
\end{table}
The observations made in the case of algebras with one product are, mutatis mutandis, also valid here, namely the
relations are trivially satisfied when one takes at least two elements in $A_1$, and the relation which is obtained
by taking all elements in $A_0$ is a consequence of the relations which are obtained by taking two elements in
$A_0$ and taking the other element in $A_1$. To see the latter claim, it suffices to consider, as in
(\ref{eq:sum_of_all}), the following formulas, which follow easily from Table \ref{tab:odot_gb},
\begin{align*}
   &(\a0\odot \b0)\odot \c1+(\a0\odot \b1)\odot \c0+(\a1\odot \b0)\odot \c0 =\ul{(a\cdot b)\cdot c}0\;,\\
   &\gb{\a0,\b0}\odot \c1+\gb{\a0,\b1}\odot \c0+\gb{\a1,\b0}\odot \c0 =\ul{\lb{a,b}\cdot c}1\;,\\
   &\gb{\a0\odot \b0,\c1}+\gb{\a0\odot \b1,\c0}+\gb{\a1\odot \b0,\c0} =\ul{\lb{a\cdot b,c}}1\;,\\
   &\gb{\gb{\a0,\b0},\c1}+\gb{\gb{\a0,\b1},\c0}+\gb{\gb{\a1,\b0},\c0}=\ul{\lb{\lb{a,b},c}}1\;,
\end{align*}
together with the four formulas, corresponding to the other parenthesizing. We have used (\ref{eq:vert_pol}) to
write the above formulas in a compact form.

\begin{example}\label{exa:poisson_dendri}
We return to the example of a Poisson algebra (see Example~\ref{exa:poisson}).  We show how to obtain the relations
which an algebra $(A,*,\circ)$ must satisfy in order to belong to the corresponding dendriform category, which we
denote by $\catPdendpol$.  We have three trilinear relations defining a Poisson algebra, namely the associativity
of ``$\cdot$'', the biderivation property and the Jacobi identity. We start with associativity of $\odot$, taking
first $\a0,\b0\in A_0$ and $\c1\in A_1$, from which we find
\begin{equation*}
  \ul{(a*b+b*a)*c}1=(\a0\odot \b0)\odot \c1=\a0\odot(\b0\odot \c1)=\ul{a*(b*c)}1\;,
\end{equation*}%
so that
\begin{equation}\label{eq:Zinbiel}
  a*(b*c)=(a*b+b*a)*c\;,
\end{equation}
for all $a,b,c\in A$, which means that $(A,*)$ is a Zinbiel algebra (see Example~\ref{exa:a_com}). Similarly,
taking $\a0,\c0\in A_0$ and $\b1\in A_1$, we find
\begin{equation*}
  \ul{c*(a*b)}1=(\a0\odot \b1)\odot \c0=\a0\odot(\b1\odot \c0)=\ul{a*(c*b)}1\;,
\end{equation*}%
so that $c*(a*b)=a*(c*b)$ for all $a,b,c\in A$, which means that $(A,*)$ is a NAP algebra (see
Example~\ref{exa:a_com}). Since every Zinbiel algebra is a NAP algebra, we don't need to state the NAP condition
for $*$. By symmetry (recall that ``$\cdot$''  is commutative) we also don't need to consider the case of $\b0,\c0\in
A_0$ and $\a1\in A_1$. Similarly, the derivation property $\lb{a\cdot b,c}=\lb{a,c}\cdot b+a\cdot \lb{b,c}$ is
symmetric in $a$ and $b$, so we get by the above procedure only two equations, which can be written in the
following symmetric form:
\begin{align}
  (a*b+b*a)\circ c&= a*(b\circ c)+b*(a\circ c)\;,\label{eq:poisson_2}\\
  (a\circ b-b\circ a)*c&=a*(b\circ c)-b\circ(a*c) \;. \label{eq:poisson_22}
\end{align}
Finally, because the Jacobi identity is symmetric in all of its variables, we get only one equation from the Jacobi
identity, namely  the pre-Lie condition
\begin{equation}\label{eq:poisson_3}
  (a\circ b-b\circ a)\circ c=a\circ(b\circ c)-b\circ(a\circ c)\;.
\end{equation}%
It follows that equations (\ref{eq:Zinbiel}) -- (\ref{eq:poisson_3}) are the four
relations of~$\catPdendpol$.

An algebra $(A, *, \circ )$ which satisfies (\ref{eq:Zinbiel}) -- (\ref{eq:poisson_3}) (i.e., an algebra
in~$\catPdendpol$) is exactly what M. Aguiar in \cite{aguiar} calls a \emph{pre-Poisson algebra}. Thus, our general
procedure to obtain $\catdendpol$ from $\catpol$ yields a canonical way to obtain the concept of a pre-Poisson
algebra from the concept of a Poisson algebra.
\end{example}

\begin{remark}
The relations which every polarized $\cC$-tridendriform algebra must satisfy are similarly obtained when the
relations $\cR_i$ are multilinear, but as in the case of $\cC$-dendriform algebras, the relations obtained by
substituting general elements from the union of $A_0$ and $A_1$ in~$\cR_i$ are all non-trivial, so there are many
more relations for a polarized $\cC$-tridendriform algebra than for a polarized $\cC$-dendriform algebra. The only
relation which we don't need to consider is the one obtained by substituting only elements from $A_0$ in $\cR_i$,
since the obtained relation is the sum of all the other relations obtained by substituting elements from the union
of $A_0$ and $A_1$ in~$\cR_i$.
\end{remark}

\begin{example}\label{exa:post-Poisson}
We continue example \ref{exa:poisson_dendri} and give the relations which an algebra $(A,*,\circ,\mid\,,\square)$
must satisfy in order to be a polarized $\cP$-tridendriform algebra.  We get the following three equations from
associativity, where the first one is obtained using the same substitutions as (\ref{eq:Zinbiel}), while the two
other equations are obtained respectively by substituting in the associativity relation two or three elements from
$A_0$:
\begin{align*}
  a*(b*c)&=(a*b)*c+(b*a)*c+(a\mid b)*c\;,\\
  a*(b\mid c)&=(a*b)\mid c\;,\\
  a\mid(b\mid c)&=(a\mid b)\mid c\;. 
\end{align*}
By symmetry, the Jacobi identity implies that we only get three relations from it, by substituting respectively
one, two or three elements from $A_0$:
\begin{align*}
  a\circ(b\circ c)-b\circ(a\circ c)&=(a\circ b-b\circ a+a\squ b)\circ c\;,\\
  (a\circ b)\squ c &=a\circ(b\squ c)+(a\circ c)\squ b\;,\\
  0&=(a\squ b)\squ c+b\squ(c\squ a)+c\squ(a\squ b)\;.
\end{align*}
Finally, the derivation property leads to the following five relations:
\begin{align*}
  a*(b\circ c)+b*(a\circ c)&=(a*b+b*a+a\mid b)\circ c\;,\\
  a*(b\circ c)-b\circ(a*c)&=(a\circ b-b\circ a+a\squ b)*c\;,\\
  (a*b)\squ c&=a*(b\squ c)+b\mid(a\circ c)\;,\\
  c\circ(a\mid b)&=a\mid(c\circ b)+b\mid(c\circ a)\;,\\
  (a\mid b)\squ c&=a\mid(b\squ c)+b\mid(a\squ c)\;.
\end{align*}
These $11$ equations are, together with the commutativity and anticommutativity of $\vert$ and $\square$, precisely
the 13 relations \cite[Eqs.\ 48--60]{poisson_bi} which define the notion of a \emph{post-Poisson algebra}.
\end{example}

\subsection{Polarization}\label{par:polarization}
We show in this subsection how the two notions of dendriform algebras, introduced in Sections \ref{sec:C-dendri}
and \ref{par:pol_dendri}, are related via a process of polarization. We first recall from \cite{MarklRemm} the
notion of polarization for an algebra $(A,\mu)$. Two new products ``$\cdot$'' and $\LB$ are defined on $A$ by
setting
\begin{equation}
  a\cdot b:=\frac12(ab+ba)\;,\qquad\hbox{and}\qquad \lb{a,b}:=\frac12(ab-ba)\;,
\end{equation}
for all $a,b\in A$ (recall that 2 is assumed invertible in the base ring $R$). This procedure is called
\emph{polarization}. Notice that we can easily reconstruct $\mu$ from the two products ``$\cdot$'' and $\LB$,
because $ab=a\cdot b+\lb{a,b}$, for all $a,b\in A$; this is what is called \emph{depolarization}. Thus, we have a
natural way to associate to each algebra $(A,\mu)$ a \emph{polarized} algebra $(A,\cdot,\LB)$ and
vice-versa. Obviously, a commutative algebra corresponds to a polarized algebra with $\LB=0$ and vice-versa, and
similarly for an anticommutative algebra, so we will only be interested in polarized algebras for which both
products are non-trivial.

\begin{example}\label{exa:P_to_Poisson}
The P-algebras introduced in Example \ref{exa:P-algebra} correspond by polarization/depolarization to Poisson
algebras, see \cite{MarklRemm}.
\end{example}

Let $\cC$ be the category of all algebras (with one product) which satisfy a given collection of relations
$\cR_1=0,\dots,\cR_k=0$. Applying polarization to all objects of $\cC$ leads to a category $\catpol$ of polarized
algebras; the morphisms in this new category are the algebra homomorphisms. Thus, by definition,
$(A,\cdot,\LB)\in\catpol$ if and only if $(A,\mu)\in\cC$, with $\mu(a,b):=a\cdot b+\lb{a,b}$ for all $a,b\in
A$. Alternatively, we can polarize the relations $\cR_i=0$ of $\cC$ by substituting in $\cR$ for $ab=\mu(a,b)$ the
sum $a\cdot b+\lb{a,b}$. Then $\catpol$ can also be described as the category of all polarized algebras, satisfying
these relations. Notice that the relations in $\cC$ are multilinear if and only if the polarized relations are
multilinear. The above polarization and depolarization procedures define inverse functors $\cC\to\catpol$ and
$\catpol\to\cC$ which make $\cC$ and $\catpol$ into isomorphic categories.

For given relations $\cR_1=0,\dots,\cR_k=0$ (in one operation) we have constructed four categories
$\cC,\,\catpol,\,\catdend$ and $\catdendpol$ and three functors, as in the following diagram, which we completed
into a square by adding a pair of inverse arrows between $\catdend$ and $\catdendpol$; the commutativity of the
diagram is easily established.
\begin{equation}\label{eq:com_pol_diagram}
  \begin{tikzcd}[row sep=8ex, column sep=15ex]
    \cC,\mu\arrow[shift left]{r}{(ab+ba)/2,(ab-ba)/2}&
      \catpol,\cdot,\LB\arrow[shift left]{l}{a\cdot b+\lb{a,b}}\\
      \catdend,\prec,\succ\arrow{u}{a\prec b+a\succ b}
                          \arrow[shift left]{r}{\frac{a\succ b+b\prec a}2,\frac{a\succ b-b\prec a}2}
    &\catdendpol,*,\circ\arrow[swap]{u}{a*b+b*a,a\circ b-b\circ a}
    \arrow[shift left]{l}{b*a-b\circ a,a*b+a\circ b}
  \end{tikzcd}
\end{equation}
In analogy with the upper arrows, we call the lower arrows \emph{polarization}
and \emph{depolarization}.
These arrows define functors which are isomorphisms of categories, just like the upper arrows. Notice that by
commutativity of the diagram, a polarized $\cC$-dendriform algebra can also be defined as an algebra $(A,*,\circ)$
whose depolarized form $(A,\prec,\succ)$ is a $\cC$-dendriform algebra (which justifies the terminology). Indeed,
according to the definition and by depolarization, $(A,*,\circ)\in\catdendpol$ if and only if $(A\times
A,\bullet)\in\cC$, with
\begin{align*}
  (a,x)\bullet (b,y)&=(a,x)\odot (b,y)+\gb{(a,x),(b,y)}\\
  &=(b*a-b\circ a+a*b+a\circ b,a*y+a\circ y+b*x-b\circ x)\\
                   &=(a\prec b+a\succ b, a\succ y+x\prec b)\;.
\end{align*}%
We have obtained exactly the condition that the depolarized form $(A,\prec,\succ)$ of $(A,*,\circ)$ belongs to
$\catdend$ (see Definition \ref{def:dendri_gen}), showing our claim.

\begin{remark}
Polarization and depolarization can also be defined for tridendriform and polarized tridendriform algebras, leading
for any category of algebras $\cC$ as above, to an isomorphism of the category $\cattrid$ of $\cC$-tridendriform
algebras and the category $\cattridendpol$ of polarized $\cC$-tridendriform algebras. On objects, the pair of
inverse isomorphisms is given by
\begin{equation}\label{eq:com_pol_diagram_2}
  \begin{tikzcd}[row sep=8ex, column sep=35ex]
      \cattrid,\prec,\succ,.\arrow[shift left]{r}{\frac{a\succ b+b\prec a}2,\frac{a\succ b-b\prec
          a}2,\frac{ab+ba}2,\frac{ab-ba}2} 
    &\cattridendpol,*,\circ,\vert,\square
    \arrow[shift left]{l}{b*a-b\circ a,a*b+a\circ b,a\mid b+a\square b}
  \end{tikzcd}.
\end{equation}
They extend the pair of lower arrows in (\ref{eq:com_pol_diagram}) and lead to a commutative diagram, as in
(\ref{eq:com_pol_diagram}).
\end{remark}
\begin{example}
We return once more to the case of P-algebras and Poisson algebras which, as we recall, correspond under
polarization; this is why we also refer to Poisson algebras as polarized P-algebras, and similarly for their
dendriform and tridendriform algebras. Specialized to this case, the above results can be summarized in the
following commutative diagram, in which the horizontal arrows are given by the horizontal arrows in
(\ref{eq:com_pol_diagram}) and~(\ref{eq:com_pol_diagram_2}):
\begin{equation*}
  \begin{tikzcd}[row sep=8ex, column sep=15ex]
    \cP,\mu\arrow[shift left]{r}&
      \catPpol,\cdot,\LB\arrow[shift left]{l}\\
    \catPtrid,\prec,\succ,.\arrow[shift left]{r}{}\arrow{u}{a\prec b+a\succ b+a.b}& 
      \catPtridendpol,*,\circ,\vert,\square\arrow[shift left]{l}{}\arrow[swap]{u}{a*b+b*a+a\vert b,a\circ b-b\circ
        a+a\squ b}\\
      \catPdend,\prec,\succ\arrow{u}{a\prec b,a\succ b,0} \arrow[shift left]{r}{}
    &\catPdendpol,*,\circ\arrow[swap]{u}{a*b,a\circ b,0,0}
    \arrow[shift left]{l}{}
  \end{tikzcd}
\end{equation*}
It was already pointed out by M. Aguiar in \cite{aguiar} that, if $(A, *, \circ )\in\catPdendpol$,
i.e., is a pre-Poisson algebra, and we define new operations on $A$ by $a\cdot b=a*b+b*a$ and $\{a, b\}=a\circ
b-b\circ a$, for all $a, b\in A$, then $(A, \cdot, \PB)$ is a Poisson algebra. It corresponds to the composition of
the two right arrows in the diagram.
\end{example}

\subsection{Application I: deformations of dendriform algebras}\label{par:defo}
In \cite{aguiar}, M. Aguiar introduced the notion of deformation for a commutative Loday dendriform algebra
$(A,\prec,\succ)$ and he showed that such a deformation makes $(A,\times,\mcirc)$ into a pre-Poisson algebra, where
$\times$ stands for $\succ$ and where the product $\mcirc$ on $A$ is constructed from the first order deformation
terms of the products $\prec$ and $\succ$. In this section we generalize this result to arbitrary $\cC$-dendriform
algebras, giving a conceptual proof of Aguiar's result.

As before, $\cC$ denotes in this section the category of all $R$-algebras satisfying a fixed set of relations
$\cR_1=0,\dots,\cR_k=0$. Let $\h$ be an indeterminate and let $\Rh$ denote the ring of formal power series
$R[[\nu]]$. More generally, for any $R$-module $A$ we denote by $\Ah$ the $\Rh$-module of formal power series
in~$\nu$ with coefficients in $A$. For a formal power series $X\in\Ah$ its evaluation at~$0$, which is the constant
term of $X$, is denoted by $X_0$.
\begin{defn}\label{def:defo}
  Let $(A,\prec_0,\succ_0)$ be a commutative $\cC$-dendriform algebra and denote $a\times b:=a\succ_0 b=b\prec_0a$
  for all $a,b\in A$. An $\Rh$-algebra $(\Ah,\prec,\succ)$ is said to be a \emph{formal deformation} of
  $(A,\prec_0,\succ_0)$ if $(\Ah,\prec,\succ)$ is a $\cC$-dendriform algebra over $\Rh$ and for any $a,b\in
  A$,
  \begin{equation*}
    (a\succ b)_0=a\succ_0b \qquad\hbox{and}\qquad
    (a\prec b)_0=a\prec_0b\;.
  \end{equation*}
  We can then define a new product on $A$ by setting, for all $a,b\in A$,
  \begin{equation}\label{equ:def_defo_term}
    a\mcirc b:=\frac{a\succ b-b\prec a}{2\h}{\Big\vert_{\h=0}}\;.
  \end{equation}%
  The algebra $(A,\times,\mcirc)$ is called the \emph{infinitesimal algebra} of the deformation. 
\end{defn}
The question which we study here is to which category the infinitesimal algebra $(A,\times,\mcirc)$ belongs. When
$\cC$ is the category of associative algebras the answer is provided by Aguiar \cite{aguiar}, who showed that
$(A,\times,\mcirc)$ is a pre-Poisson algebra. 

In order to answer the above question in general, we first introduce a few more notions and notations. Let $M$ be a
monomial which involves the (commutative and anticommutative) products ``$\cdot$'' and $\LB$ only. We define the
\emph{weight} of $M$ as the number of operations $\LB$ in $M$. Similarly, for a monomial $M$ in the products $*$
and $\circ$, its weight is the number of operations $\circ$ in~$M$. In either case, a sum $\cR$ of such monomials
is said to be \emph{homogeneous of weight $m$} if each of its terms has weight $m$. The lowest weight part of $\cR$
is denoted by~$\underline\cR$. Finally, we denote by $\ucatpol$ (resp.\ by $\ucatdendpol$) the category of all
$R$-algebras satisfying all relations $\uR=0$, where $\cR$ runs through the linear space of relations of $\catpol$
(resp.\ of $\catdendpol$).
\begin{prop}\label{prp:defo_homog}
  Let $(\Ah,\prec,\succ)$ be a formal deformation of a commutative algebra $(A,\prec_0,\succ_0)\in\catdend$, with
  deformation algebra $(A,\times,\mcirc)$. Then
  \begin{equation}\label{eq:defo_alg}
    (A,\times,\mcirc)\in\ucatdendpol\;.
  \end{equation}
  In particular, when the relations of $\catdendpol$ are generated by weight homogeneous relations, then
  $(A,\times,\mcirc)\in\catdendpol$. Also, when the relations of $\cC$ are multilinear,
  $\ucatdendpol=\(\ucatpol\)^{\hbox{\tiny dend}}$, so that $(A,\times,\mcirc)=\(\ucatpol\)^{\hbox{\tiny dend}}$.
\end{prop}
\begin{proof}
We will only prove here that $(A,\times,\mcirc)\in\ucatdendpol$ leaving the more technical proof that
$\ucatdendpol=\(\ucatpol\)^{\hbox{\tiny dend}}$ to the end of the section.

Given a formal deformation $(\Ah,\prec,\succ)$ we can construct by polarization (which, as we recall, is an
isomorphism of categories) an algebra $(\Ah,*,\circ)$, which is a polarized dendriform algebra over $\Rh$. We
define new products $*_i$ and $\circ_i$ on $A$ by setting for all $a,b\in A$,
\begin{align}\label{equ:defo_pol}
  a*b&=a*_0b+a*_1b\,\h+a*_2b\,\h^2+\cdots\;,\nonumber\\
  a\circ b&=a\circ_0b+a\circ_1b\,\h+a\circ_2b\,\h^2+\cdots\;.
\end{align}
Since, by polarization, $a\circ b=(a\succ b-b\prec a)/2$ and $a* b=(a\succ b+b\prec a)/2$ (see
(\ref{eq:com_pol_diagram})), we have by commutativity of $(A,\prec,\succ)$ that $a*_0b=a\times b$ and that
$a\circ_0b=0$; also, the definition of $\mcirc$ implies that $a\circ_1b=a\mcirc b$ for all $a,b\in A$. Hence,
(\ref{equ:defo_pol}) can be rewritten as
\begin{align}
  a*b&=a\times b+a*_1b\,\h+a*_2b\,\h^2+\cdots\;,\label{equ:defo_pol_2}\\
  a\circ b&=\hskip 1.3cm a\mcirc b\,\h+a\circ_2b\,\h^2+\cdots\;,\label{equ:defo_pol_3}
\end{align}
where the dots stand for terms containing $\h^i$ with $i>2$. Suppose now that $\cR=0$ is a relation of
$\catdendpol$. Writing $\cR$ as $\cR_{*,\circ}$ to indicate the products which are involved, we may also consider
$\cR_{\times,\mcirc}$. We need to show that $\uR_{\times,\mcirc}(a_1,\dots,a_n)=0$ for all $a_1,\dots,a_n\in A$.
To do this, consider the relation $\cR_{*,\circ}(a_1,\dots,a_n)=0$. In view of (\ref{equ:defo_pol_2}) and
(\ref{equ:defo_pol_3}),
\begin{equation}\label{equ:R_to_R}
  \cR_{*,\circ}(a_1,a_2,\dots,a_n)=\uR_{\times,\mcirc}(a_1,a_2,\dots,a_n)\h^d+\cdots,
\end{equation}
where $d$ denotes the lowest weight of the terms of $\cR$, i.e.\ the weight of~$\uR$. It follows that
$(A,\times,\mcirc)$ satisfies the relation $\uR_{\times,\mcirc}=0$, as was to be shown.
\end{proof}
\begin{example}\label{exa:defo_assoc}
Let $\cC$ be the category of all associative algebras (over $R$). Then, by polarization, the following are the
relations in $\catpol$ (see \cite{MarklRemm}):
\begin{align}
  \lb{a\cdot b,c}&=a\cdot\lb{b,c}+\lb{a,c}\cdot b\;,\label{equ:cat_ass_1}\\  
  \lb{\lb{a, b},c}&=(b\cdot c)\cdot a-(c\cdot a)\cdot b\;.\label{equ:cat_ass_2}
\end{align}
Recall that (\ref{equ:cat_ass_2}) implies the Jacobi identity, which is weight homogeneous (of weight 2), just
like the derivation property (\ref{equ:cat_ass_1}) (of weight~1). Notice that the lowest weight part of
(\ref{equ:cat_ass_2}) is $(b\cdot c)\cdot a=(c\cdot a)\cdot b$, which is associativity (since ``$\cdot$''
commutative). It follows that $\ucatpol$ is the category of Poisson algebras, hence that $\ucatdendpol$ is the
category of pre-Poisson algebras. This shows that the infinitesimal algebra of a deformation of a Loday
dendriform algebra is a pre-Poisson algebra, as was first shown by Aguiar~\cite{aguiar}.
\end{example}

\begin{example}\label{exa:homog}
The relations which define Poisson algebras (see Example~\ref{exa:poisson}) are 3-linear and homogeneous:
associativity is of weight 0, the derivation property is of weight 1 and the Jacobi identity is of weight two. For
$A_3$-associative algebras and $LA$-algebras in polarized form, the relations are also easily written in
homogeneous form.  It follows that the second part of Proposition \ref{prp:defo_homog} can be applied to these
algebras: in each of these cases, the infinitesimal algebra $(A,\times,\mcirc)$ of the deformation belongs to
$\catdendpol$.
\end{example}

\begin{remark}
Proposition \ref{prp:defo_homog} is easily adapted to the classical case of formal deformations $(A,\mu)$ of
commutative algebras $(A,\mu_0)\in\cC$. The infinitesimal algebra is then defined as $(A,\mu_0,\diamond)$, where
\begin{equation*}
  a\diamond b:=\frac{\mu(a,b)-\mu(b,a)}{2\h}{\Big\vert_{\h=0}}\;.
\end{equation*}%
One shows as in the proof of Proposition \ref{prp:defo_homog} that $(A,\mu_0,\diamond)\in\ucatpol$. In the case of
associative algebras, $\ucatpol$ is the category of Poisson algebras (see Example \ref{exa:defo_assoc}), so we
recover the classical result that the infinitesimal algebra of a deformation of an associative algebra is a Poisson
algebra.
\end{remark}

\begin{remark}
One may also consider more generally deformations of $\cC$-tridendriform algebras. Recall that in a commutative
$\cC$-tridendriform algebra $(A,\prec,\succ,.)$, one also requires the last product to be commutative. The weight
of a relation $\cR=\cR_{*,\circ,\vert,\squ}$ is now defined such that $*$ and $\vert$ have weight 0, while $\circ$
and $\squ$ have weight 1. It is clear that all the above results generalize to this case. The infinitesimal algebra
has now four operations. For example, when $\cC$ is the category of associative algebras, the infinitesimal algebra
is a post-Poisson algebra (see Example \ref{exa:post-Poisson}).
\end{remark}

\begin{remark}
We have considered deformations of commutative dendriform algebras, but everything can be easily adapted to
anticommutative dendriform algebras: the rôles of $*$ and $\circ$ are exchanged in the sense that one will have now
that $*_0=0$, that $*_1=\times$ and $\circ_0=\mcirc$, where $(A,\mcirc)$ is the original anticommutative dendriform
algebra (written as an algebra with one operation). As we have seen in Section \ref{par:comm_anticom},
$A_3$-associative, $LA$ and P-algebras which are anticommutative are Lie algebras, so there are many natural examples
of this case.
\end{remark}

To finish this section, we prove that when the relations of $\cC$ are multilinear,
$\ucatdendpol=\(\ucatpol\)^{\hbox{\tiny dend}}$, as stated in (\ref{eq:defo_alg}). The property says that the
lowest weight parts of all relations in $\ucatdendpol$ are obtained by dendrifying the lowest weight parts of all
relations in $\ucatpol$. Notice that since each dendrification of a monomial of weight $k$ (involving the products
``$\cdot$'' and $\LB$ only) is homogeneous of weight $k$, one has that all algebras in $\(\ucatpol\)^{\hbox{\tiny
dend}}$ are also algebras of $\ucatdendpol$. We therefore only need to prove the reciprocal inclusion.

Notice also that we may restrict ourselves to $n$-linear relations, for a fixed $n$, since the
dendrification of a $k$-linear relation is $k$-linear, i.e. we may suppose that all relations $\cR_1,\dots,\cR_k$
of $\catpol$, and hence also of $\catdendpol$ are $n$-linear.

For $0\leqslant\ell\leqslant n$, consider the free $R$-modules $\cM_\ell$ and $\tcM_\ell$, generated by all
$\ell$-linear monomials $M$ involving only the (commutative and anticommutative) products ``$\cdot$'' and $\LB$
only, respectively generated by all $\ell$-linear monomials $\tM$ involving only the products $*$ and $\circ$ in
$n$ variables, say $x_1,\dots,x_n$. Their direct sums are denoted $\cM$ and $\tcM$ respectively. Elements of
$\cM_\ell$ and $\tcM_\ell$ are also said to be of \emph{length} $\ell$; notice that the weight of a monomial of
length $\ell$ is between $0$ and $\ell-1$ (included). The modules ${\cM}_\ell$ and $\tcM_\ell$ admit natural
decompositions
\begin{equation*}
  {\cM}_\ell={\cM}^0_\ell\oplus\cdots\oplus{\cM}^{\ell-1}_\ell\quad\hbox{and}\quad
  {\tilde\cM}_\ell={\tilde\cM}^0_\ell\oplus\cdots\oplus{\tilde\cM}^{\ell-1}_\ell\;,
\end{equation*}
where ${\cM}^i_\ell\subset \cM_\ell$ and ${\tcM}^i_\ell\subset\tcM_\ell$, are the submodules generated by the
monomials of weight~$i$.  Each monomial $M$ of ${\cM}_\ell$ of length at least two can be decomposed as $M=M_1\cdot
M_2$ or $M=\lb{M_1,M_2}$; this decomposition is unique up to the order of the factors.


We describe the process of dendrification of multilinear relations of $\catpol$, introduced and studied in Section
\ref{sec:multilinear_pol}, in terms of the linear maps
\begin{equation*}
  \varphi_0,\varphi_1,\dots,\varphi_n:\cM\to\tcM\;,
\end{equation*}
which we define on monomials $M$, using induction on the length of $M$: 
\begin{equation*}
  \varphi_0(M):=
    \left\{
    \begin{array}{ll}
      x_i&\text{ if }M=x_i\;;\\
      \varphi_0(M_1)*\varphi_0(M_2)+\varphi_0(M_2)*\varphi_0(M_1)&\text{ if }M=M_1\cdot M_2\;;\\
      \varphi_0(M_1)\circ\varphi_0(M_2)-\varphi_0(M_2)\circ\varphi_0(M_1)&\text{ if }M=\lb{M_1,M_2}\;,\end{array}
    \right.
\end{equation*}
and for $p=1,\dots,n$ we define
\begin{equation*}
  \varphi_p(M):=
    \left\{
      \begin{array}{ll}
        0\quad& \text{if $M$ is independent of $x_p$\;;}\\
        x_p&\text{if $M=x_p$}\\
        \varphi_0(M_1)*\varphi_p(M_2)&\text{if $M=M_1\cdot M_2$ and $M_2$ depends on $x_p$\;;}\\
        \varphi_0(M_1)\circ\varphi_p(M_2)&\text{if $M=\lb{M_1,M_2}$ and $M_2$ depends on $x_p$\;.}
      \end{array}
    \right.
\end{equation*}
It is clear that these maps are well-defined and that they preserve the length and the weight of a monomial. Notice
that, by construction, in all terms of $\varphi_p(\cM)$ the variable $x_p$ is located at the last
position. Therefore, the images of the maps $\varphi_1,\dots,\varphi_n$ are in direct sum.

To see the relation with dendrification, let $\cR=0$ be an $n$-linear relation of~$\catpol$. Then $\cR\in\cM$ and
for $p=1,\dots,n$, the relation $\varphi_p(\cR)=0$ is precisely the relation obtained by substituting in
$\cR_{\cdot,\LB}$ for the $p$-th variable $(0,x_p)$ and for the $q$-th variable $(x_q,0)$, where $q\neq p$.

\begin{lemma}\label{lma:injectivity}
  The maps $\varphi_0,\dots,\varphi_n$ are injective.
\end{lemma}

\begin{proof}
Let $\tM$ be a monomial of $\tcM$. We show that there exists a unique monomial
$M\in\cM$ such that $\tM$ is a term of $\varphi_0(M)$; from it the injectivity of $\varphi_0$ is clear.

We do this by induction on the length of $\tM$. When $\tM$ is of length 1, the claim is trivially true, so let us
assume that the claim is true for monomials of length strictly less than some $\ell\geqs2$. Let $\tM$ be a monomial
of $\tcM$ of length $\ell$. We can write $\tM$ uniquely as $\tM=\tM_1*\tM_2$ or $\tM=\tM_1\circ\tM_2$, up to the
order of the factors. By the induction hypothesis there exists a unique couple $(M_1,M_2)$ such that $\tM_1$ and
$\tM_2$ are terms of $\varphi_0(M_1)$ and $\varphi_0(M_2)$ respectively, and hence such that $\tM$ is a term of
$\varphi_0(M_1)*\varphi_0(M_2)$ or $\varphi_0(M_1)\circ\varphi_0(M_2)$, depending on whether $\tM=\tM_1*\tM_2$ or
$\tM=\tM_1\circ\tM_2$. It follows that, if we define $M:=M_1\cdot M_2$ or $M:=\lb{M_1,M_2}$, depending on whether
$\tM=\tM_1*\tM_2$ or $\tM=\tM_1\circ\tM_2$, then $\varphi_0(M)=\tM$. Since the decomposition of $\tM$ is unique up
to the order of the factors, $M$ is unique. This shows the claim, and hence the injectivity of $\varphi_0$.

In order to show the injectivity of the other maps $\varphi_1,\dots,\varphi_n$ one proceeds in a similar way: one
shows as above that given any monomial $\tM$ of $\tcM$ there exists a unique monomial $M$ of $\cM$ and a unique
integer $p\in\set{1,\dots,n}$ such that $\tM$ is a term of $\varphi_p(M)$.
\end{proof}

\begin{lemma}\label{prp:homog_iff}
Let $\cR_1,\dots,\cR_k\in\cM_n$. For any constants $\l_i^p\in R$ $(1\leqslant i\leqs k$ and $p=1,\dots,n)$, not all
equal to zero,
\begin{equation}\label{eq:lowest_weight}
  \underline{\sum\nolimits_{i=1}^k\sum\nolimits_{p=1}^n\lambda_i^p \varphi_p(\cR_i)}=
  \sum_{p=1}^n\varphi_p\left(\underline{\sum\nolimits_{i=1}^k\lambda_i^p\cR_i}\right)\;.
\end{equation}
\end{lemma}

\begin{proof}
For $i=1,\dots,k$, let $\cR_i=\cR^0_i+\cdots+\cR^{n-1}_i$ be the weight decomposition of $\cR_i$.  By $R$-linearity
of the maps $\varphi_p$,
\begin{equation*}
  \sum_{i=1}^k\sum_{p=1}^n\lambda_i^p \varphi_p(\cR_i)=
  \sum_{\ell=m}^{n-1}A_\ell\;,
  \qquad
  \text{where}
  \qquad
  A_\ell=\sum_{p=1}^n\varphi_p(\sum_{i=1}^k\lambda_i^p \cR^\ell_i)\;,
\end{equation*}
and where $m$ is chosen such that $A_0,\dots,A_{m-1}=0$ and $A_{m}\neq0$.  Since the maps $\varphi_p$ are
weight-preserving, $A_\ell$ is homogeneous of weight $\ell$, and so $A_m$ is equal to the left hand side of
(\ref{eq:lowest_weight}). Let $0\leqslant\ell<m$. Then $\sum_{p=1}^n\varphi_p(\sum_{i=1}^k\lambda_i^p
\cR^\ell_i)=A_\ell=0$, so that $\varphi_p(\sum_{i=1}^k\lambda_i^p \cR^\ell_i)$ for all $p$, since the images of the
maps $\varphi_1,\dots,\varphi_n$ are in direct sum. Since the maps $\varphi_p$ are injective (Lemma
\ref{lma:injectivity}), this implies that $\sum_{i=1}^k\lambda_i^p \cR^\ell_i=0$ for $\ell=0,\dots,m-1$.  Also,
$\sum_i\lambda_i^p \cR^m_i\neq0$ since $A_m\neq0$. It follows that
\begin{equation*}
  \underline{\sum\nolimits_{i=1}^k\lambda_i^p\cR_i}=\underline{\sum\nolimits_{i=1}^k\sum\nolimits_{\ell=0}^{n-1}
    \lambda_i^p\cR_i^\ell}  ={\sum\nolimits_{i=1}^k\lambda_i^p\cR_i^m}\;,
\end{equation*}
so that $A_m$ is also equal to the right hand side of (\ref{eq:lowest_weight}). 
\end{proof}
We use Lemma \ref{prp:homog_iff} to show that all algebras in $\ucatdendpol$ are also algebras of
$\(\ucatpol\)^{\hbox{\tiny dend}}$, so that $\ucatdendpol=\(\ucatpol\)^{\hbox{\tiny dend}}$. Suppose that
$\cR_1=0,\dots,\cR_k=0$ is a basis for the module of all $n$-linear relations of $\catpol$.
Let $\cR=0$ be a relation of $\ucatdendpol$. By definition, $\cR$ is the lowest weight part of
$\sum\nolimits_{i=1}^k\sum\nolimits_{p=1}^n\lambda_i^p \varphi_p(\cR_i)$, for some constants $\lambda_i^p$.  In
view of the lemma, $\cR$ is obtained by dendrification of some relations in $\ucatpol$, namely the $p$ relations
$\sum_{i=1}^k\lambda_i^p\cR_i=0$, for $p=1,\dots,n$. This shows that $\cR=0$ is a relation of
$\(\ucatpol\)^{\hbox{\tiny dend}}$.

\subsection{Application II: filtered dendriform algebras}\label{par:filtrations}
As a second application of polarized dendriform algebras, we generalize another result of Aguiar \cite{aguiar},
which is itself an analogue for Loday dendriform algebras of the well-known result which says that the graded
algebra associated to an almost commutative filtered associative algebra is a Poisson algebra.

Let $(A,\prec,\succ)$ be an algebra. An (increasing) (\emph{filtration} on $A$ is an increasing sequence of
subspaces $A_0\subseteq A_1\subseteq A_2 \subseteq \cdots$ such that
$$
  A=\bigcup_{i\geqs 0}A_i \qquad \hbox{and} \qquad (A_i\prec A_j+A_i\succ A_j)\subseteq A_{i+j}, 
$$
for all $i,j\geqs0$. Then $A$ is called a \emph{filtered algebra}. It is convenient to set $A_{i}:=\set0$
for $i<0$. The associated graded algebra is, as an $R$-module,
$$
  \Gr(A):=\bigoplus _{i\geqs 0}\frac{A_{i}}{A_{i-1}}
$$
and inherits two products from $\prec$ and $\succ$, which are still denoted by  $\prec$ and~$\succ$. They are
(well-) defined by setting, for $a\in A_i$ and $b\in A_j$, with $i,j\geqs0$,
$$
  (a+A_{i-1})\prec (b+A_{j-1}):=(a\prec b+A_{i+j-1})\in \frac{A_{i+j}}{A_{i+j-1}}\;,
$$
and similarly for $\succ$. As in the case of algebras with one operation, $A$ and $\Gr(A)$ are canonically
isomorphic as $R$-modules, but not as algebras. It is however clear that any $n$-linear relation which is satisfied
by the original products $\prec$ and $\succ$ will be satisfied by the induced products.

We will be interested in \emph{almost commutative} filtered algebras, which have the property that the associated
graded algebra is commutative, i.e., $a\prec b=b\succ a$ for all $a,b\in\Gr(A)$. As before, we then view $\Gr(A)$
as an algebra with one operation $\times$ (setting as usual $\times:=\ \succ$), and $\Gr(A)$ can be equipped with
another product, defined for $a\in A_i$ and $b\in A_j$, with $i,j\geqs0$ by
\begin{equation}\label{equ:def_filt_term}
  (a+A_{i-1})\mcirc (b+A_{j-1}):=(a\succ b-b\prec a +A_{i+j-2})\in \frac{A_{i+j-1}}{A_{i+j-2}}\;.
\end{equation}
The question is now again to which category $(\Gr(A),\times,\mcirc)$ belongs. When $\cC$ is the category of
associative algebras, Aguiar's answer is that $(\Gr(A),\times,\mcirc)$ is a pre-Poisson algebra, as in the case of
deformations (see~\cite{aguiar}). We will give here the answer for arbitrary algebras; as we will see, the result
is very similar to the result which we obtained for deformations (Section \ref{par:defo}). The definitions and
assumptions are the same as in the latter section, except that the relations of $\cC$ (and hence of $\catpol$) are
supposed here to be multilinear.

\begin{prop}\label{prp:filtered}
  Suppose that the relations of $\cC$ are multilinear. Let $(A=\cup_iA_i,\prec,\succ)$ be a commutative filtered
  algebra in $\catdend$. On $\Gr(A)$, consider the product $\times$, defined for $a,b\in\Gr(A)$ by $a\times
  b:=a\succ b$, as well as the product $\mcirc$, defined by (\ref{equ:def_filt_term}). Then
  $$(\Gr(A),\times,\mcirc)\in\ucatdendpol=\(\ucatpol\)^{\hbox{\tiny dend}}\;.$$
\end{prop}
\begin{proof}
As in the proof of Proposition \ref{prp:defo_homog}, we use polarization to transform the deformation into an
algebra of $\catdendpol$. Namely, by polarization, we have a filtered algebra $(A,*,\circ)\in\catdendpol$, having
the property that
\begin{equation}\label{equ:inclusion}
  A_i*A_j\subset A_{i+j}\;,\qquad\hbox{and}\qquad A_i\circ A_j\subset A_{i+j-1}\;.
\end{equation}
In terms of $*$ and $\circ$, the above definitions of $\times$ and $\mcirc$ now amount to setting, for $a\in A_i$
and $b\in A_j$,
\begin{align}
  (a+A_{i-1})\times (b+A_{j-1})&:=a*b+A_{i+j-1}\;,\label{eqn:times}\\
  (a+A_{i-1})\mcirc (b+A_{j-1})&:=a\circ b+A_{i+j-2}\;.\label{eqn:mcirc}
\end{align}%
Suppose now that $\cR=\cR_{*,\circ}$ is an $n$-linear relation of $\catdendpol$ and recall that we denote the
lowest weight part of $\cR$ by $\uR$. The weight of $\uR$ is denoted by $d$. Let $a_1,a_2,\dots,a_n\in A$ with
$a_i\in A_{j_i}$ for $i=1,\dots,n$. Then
\begin{align*}
  \lefteqn{\ul{\cR}{\times,\mcirc}(a_1+A_{j_1-1},\dots,a_n+A_{j_n-1})}\\
  &=\ul{\cR}{*,\circ}(a_1,\dots,a_n)+A_{j_1+\cdots+j_n-d-1}\\
  &={\cR}_{*,\circ}(a_1,\dots,a_n)+A_{j_1+\cdots+j_n-d-1}
  =A_{j_1+\cdots+j_n-d-1}\;.
\end{align*}
where we used in the last step that $(A,*,\circ)$ satisfies $\cR$.  It follows that $(\Gr(A),\times,\mcirc)$
satisfies the relation $\uR=0$. Therefore, $(\Gr(A),\times,\mcirc)$ satisfies all relations of $\ucatdendpol$, and
so ($\Gr(A),\times,\mcirc)\in\ucatdendpol$.
\end{proof}
\begin{example}
We return once more to the case where $\cC$ is the category of associative algebras. We have already analyzed the
relations defining $\catdendpol$ in Example \ref{exa:defo_assoc} where we have shown that the lowest weight terms
of the relations are the relations which define a pre-Poisson algebra. Hence, we find that if $(A,\prec,\succ)$ is
an almost commutative filtered Loday dendriform algebra, then $(\Gr(A),\times,\mcirc)$ is a pre-Poisson algebra. We
thereby recover Aguiar's result, cited above.
\end{example}
The strong similarity between our results on filtrations and on deformations is not accidental. Indeed, let
$(\Ah,\prec,\succ)$ be a formal deformation of a commutative algebra $(A,\prec_0,\succ_0)\in\catdend$, where we
assume that the relations which define $\cC$ are multilinear. Setting $\Ah_i:=\h^i\Ah$ for all $i\in\bbN$ it is
clear that $(\Ah,\prec,\succ)$ is a filtered $\cC$-dendriform algebra. Notice that the filtration is
\emph{descending}, so that $\Gr(\Ah)$ is now defined as $\Gr(\Ah):=\bigoplus_{i\geqs0}{\Ah_i}/{\Ah_{i+1}}$, and
that $\Gr(\Ah)$ is commutative.  Though ascending and descending filtrations (indexed by $\bbN$) are from many
points of view different, it is easily verified that the above results on ascending filtrations hold also for
descending filtrations. In particular, $(\Gr(\Ah),\times,\mcirc)\in\ucatdendpol$, as in Proposition
\ref{prp:filtered}. Under the canonical isomorphisms ${\Ah_i}/{\Ah_{i+1}}\simeq A$, valid for all $i\in\bbN$, we
get that $(A,\times,\mcirc)\in\ucatdendpol$, where the latter products on $A$ are inherited from the products on
$\Gr(A)$. It is easily checked that $(A,\times,\mcirc)$ is the deformation algebra of $(\Ah,\prec,\succ)$. This
shows that under the extra assumption that the relations defining $\catdendpol$ are multilinear, Proposition
\ref{prp:defo_homog} is a consequence of \ref{prp:filtered}. It should now be clear that all remarks made in
Section \ref{par:defo} also apply to almost commutative (or anticommutative) filtered algebras (always under the
assumption that the relations defining $\catdendpol$ are multilinear).

\medskip

\begin{center}
ACKNOWLEDGEMENTS
\end{center}

Parts of this paper have been written while the second author was a Visiting Professor at the University of
Poitiers in June 2018.